\documentclass[11pt]{article}
\usepackage{amsfonts}
\usepackage{amsmath}
\usepackage{enumitem}
\usepackage{geometry}                
\usepackage{graphicx}
\usepackage{caption}
\usepackage{amssymb}
\usepackage{epstopdf}
\usepackage{hyperref}             
\title{On the analysis of mixed-index time fractional differential equation systems}
\author{Kevin Burrage\footnotemark[1] \footnotemark[2] \and Pamela M. Burrage\footnotemark[1] \footnotemark[2] \and Ian W. Turner\footnotemark[1] \footnotemark[2] \and Fanhai Zeng\footnotemark[2]}
\date{}

\begin{document}
\maketitle

\renewcommand{\thefootnote}{\fnsymbol{footnote}}

      \footnotetext[1]{ARC Centre of Excellence for Mathematical and Statistical Frontiers,
      Queensland University of Technology, Australia.}
  \footnotetext[2]{School of Mathematical Sciences, Queensland University of Technology (QUT), Australia. kevin.burrage@qut.edu.au}

   \renewcommand{\thefootnote}{\arabic{footnote}}

\abstract
{In this paper we study  the class of mixed-index time fractional differential equations in which different components of the problem have different time fractional derivatives on the left hand side. We prove a theorem on the solution of the linear system of equations, which collapses to the well-known Mittag-Leffler solution in the case the indices are the same, and also generalises the solution of the so-called linear sequential class of time fractional problems. We also investigate the asymptotic stability properties of this class of problems using Laplace transforms and  show how Laplace transforms can be used to write solutions as linear combinations of generalised Mittag-Leffler functions in some cases.   Finally we illustrate our results with some numerical simulations.}


\vspace{5mm}



\section{Introduction}
Time fractional and space fractional differential equations are increasingly used as a powerful modelling tool for understanding the role of heterogeneity in modulating function in such diverse areas as cardiac electrophysiology \cite{ref1,ref1a,ref1b}, brain dynamics \cite{ref2}, medicine \cite{ref3}, biology \cite{ref4}, \cite{ref5}, porous media \cite{ref6}, \cite{ref7} and physics \cite{ref8}. Time fractional models are typically used to model subdiffusive processes (anomalous diffusion \cite{ref9}, \cite{ref10}), while space fractional models are often associated with modelling processes occurring in complex spatially heterogeneous domains \cite{ref1}.

Time fractional models typically have solutions with heavy tails as described by the Mittag-Leffler matrix function \cite{ref11} that naturally occurs when solving time fractional linear systems. However such models are usually only described by a single fractional exponent, $\alpha$, associated with the fractional derivative. The fractional exponent can allow the coupling of different processes that may be occurring in different spatial domains by using different fractional exponents for the different regimes. One natural application here would be the coupling of models describing anomalous diffusion of proteins on the plasma membrane of the cell with the behaviour of other proteins in the cytosol of the cell. Tian et al \cite{ref12} addressed this problem by coupling a stochastic model (based on the Stochastic Simulation Algorithm \cite{refSSA}) for the plasma membrane with systems of ordinary differential equations describing reaction cascades within the cell.  It may also be necessary to couple more than two models and so in this paper we introduce a formulation that focuses on coupling an arbitrary number of domains in which dynamical processes are occurring described by different anomalous diffusive processes. This leads us to consider the $r$ index time fractional differential equation problem in Caputo form
\begin{equation}
D_t^{\alpha_i} y_i = \sum_{j=1}^r A_{ij} y_j + F_i(y), \quad y_i(0) = z_i, \quad y_i \in \mathbb{R}^{m_i}, \> i=1,\cdots,r, 
 \label{eq:1}
\end{equation}
or in vector form
$$D_t^\alpha \, y = A\,y + F(y).$$
Here the $A_{ij}$ are $m_i \times m_j$  matrices, while $A$ is the associated block matrix of dimension $\sum_{j=1}^r m_j$ and $\alpha = (\alpha_1,\cdots,\alpha_r)^\top$ has all components $\alpha_i \in (0,1]$.

We believe that a modelling approach based on this formulation has not been fully developed before. We note that scalar linear sequential fractional problems have been considered whose solution can be described by multi-indexed Mittag-Leffler functions \cite{ref13}, and  there are a number of articles on the numerical solution of multi-term fractional  differential equations \cite{ref_f1,ref_f2,ref_f3}, and while mixed index problems can, in some cases, be written in the form of linear sequential problems, namely $\sum_{i=1}^R D_t^{\beta_i}\,y = f(y)$, we claim that it is inappropriate to do so in many cases.

Therefore in this paper we develop a new theorem that gives the analytical solution of equations such as (\ref{eq:1}) that reduces to the Mittag-Leffler expansion in the case that all the indices are the same (section 3) and generalises the class of linear sequential problems. We then analyse the asymptotic stability properties of these mixed index problems using Laplace transform techniques (section 4), relating our results with known results that have been developed in control theory. In section 5 we show that, in the case that the $\alpha_i$ are all rational, the solutions to the linear problem can be written as a linear combination of generalised Mittag-Leffler functions, again using ideas from control theory and transfer functions. In section 6 we present some numerical simulations illustrating the results in this paper and give some discussion on how these ideas can be used to solve semi-linear problems  either by extending the methodology of exponential integrators to Mittag-Leffler functions, or by writing the solution as sums of certain Mittag-Leffler expansions.

\section{Background}
We consider the linear system given in (\ref{eq:1}) with $r = 2$. It will be convenient to let
\begin{equation} A = \left( \begin{array}{ll} A_1 & A_2 \\ B_1 & B_2 \\ \end{array} \right), \quad y^{\top} = ( y_1^{\top}, y_2^{\top}), \quad z^{\top} = ( z_1^{\top}, z_2^{\top})
\label{eq:sec2eq1}
\end{equation}
where $A$ is $m \times m, \> m = m_1 + m_2.$ We will call such a system a time fractional index-2 system. Here the Caputo time fractional derivative with starting point at $t=0$ is defined (see Podlubny \cite{ref14}, for example), as
$$D_t^{\alpha} y(t) = \frac{1}{\Gamma(1-\alpha)} \int_0^t \frac{y'(s)}{(t-s)^{\alpha}}ds, \quad 0<\alpha < 1.$$
Furthermore, given a fixed mesh of size $h$ then a first order approximation of the Caputo derivative \cite{ref15} is given by 
$$D_t^{\alpha} y_n = \frac{1}{\Gamma(2-\alpha) h^{\alpha}} \sum_{j=1}^n (j^{1-\alpha}-(j-1)^{1-\alpha})(y_{n-j-1}-y_{n-j}).$$

If $\beta$ = $\alpha$ then the solution to (\ref{eq:1}) is given by the Mittag-Leffler expansion
\begin{equation}
y(t) = E_{\alpha}(t^{\alpha}A)\,y(0), \quad 
E_{\alpha}(z) = \sum_{j=0}^{\infty} \frac{z^j}{\Gamma(1+j \alpha)}
\label{eq:2}
\end{equation}
where $\Gamma(x)$ is the Gamma function.

If the problem is completely decoupled, say $A_2 = 0$, then from (\ref{eq:2}) the solution to (\ref{eq:1}) and (\ref{eq:sec2eq1}) satisfies
\begin{eqnarray}
y_1(t)  &=& E_{\alpha}(t^{\alpha}A_1)\, z_1 \nonumber \\
D_t^{\beta} y_2 &=& B_2 y_2 + B_1 E_{\alpha}(t^{\alpha}A_1)\, z_1.
\label{eq:4}
\end{eqnarray}
In order to solve (\ref{eq:4}), this requires us to solve problems of the form
\begin{equation}
D_t^{\beta} y_2 = B_2 y_2 + f(t).
\label{eq:5}
\end{equation}

Before making further headway, we need some additional background material.

\vspace{5mm}
\noindent
\textbf{Definition 1.}
Generalisations of the Mittag-Leffler functions are given by
\begin{eqnarray*}
E_{\alpha, \beta}(z) &=& \sum_{k=0}^{\infty} \frac{z^k}{\Gamma(\alpha k + \beta)}, \quad \texttt{Re}(\alpha) > 0 \\
E_{\alpha, \beta}^{\gamma}(z) &=& \sum_{k=0}^{\infty} \frac{(\gamma)_k}{\Gamma(\alpha k + \beta)} \, \frac{z^k}{k!}, \quad \gamma \in \mathbb{N}_0, 
\end{eqnarray*}
where $(\gamma)_k$ is the Pochhammer symbol
$$ (\gamma)_0 = 1, \quad (\gamma)_k = \gamma (\gamma+1) \cdots (\gamma+k-1).$$

\noindent
\textbf{Remark.} $E_{\alpha,1}(z) = E_{\alpha}(z), \quad E_{\alpha, \beta}^{1}(z) = E_{\alpha, \beta}(z), \quad E_1(z) = e^z.$
\vspace{5mm}

\noindent
\textbf{Lemma 1.}
$$ \left( \frac{d}{dz} \right)^n \, E_{\alpha, \beta}(z) = n! E_{\alpha, \beta+\alpha n}^{n+1} (z), \> n \in \mathbb{N}.$$

The following result will be important in section 5.

\vspace{5mm}
\noindent
\textbf{Lemma 2.} The Laplace transform of $E_{\alpha,\beta}(\lambda t^\alpha)$ satisfies
\begin{equation}
X(s) = \frac{s^\alpha}{s^\beta(s^\alpha-\lambda)}.
\label{eq:paper2eq13}
\end{equation}

\vspace{5mm}

The Caputo derivatives satisfy the following relationships.

\vspace{5mm}

\noindent
\textbf{Lemma 3.}
\begin{enumerate}
\item $ D_t^{\alpha} I^{\alpha} y(t) = y(t)$
\item $I^{\alpha} D_t^{\alpha} y(t) = y(t) - y(0)$
\item $D_t^{\alpha} y(t) = \frac{1}{\Gamma(1-\alpha)} \int_0^t \frac{y'(s)}{(t-s)^{\alpha}} ds = I^{1-\alpha}D_t y(t).$
\end{enumerate}

\vspace{5mm}

\noindent
\textbf{Lemma 4.}
The solution of the scalar, linear, non-homogeneous problem
\begin{equation}
D_t^{\alpha} y(t) = \lambda y(t) + f(t), \quad y(0) = y_0
\label{eq:6}
\end{equation}
is
\begin{equation}
y(t) = E_{\alpha}(\lambda t^{\alpha}) y_0 + \int_0^t (t-s)^{\alpha - 1} E_{\alpha \alpha} (\lambda (t-s)^{\alpha}) f(s) ds.
\label{eq:7}
\end{equation}

\noindent
\textbf{Proof:}
Using the integral form from Lemma 3, (\ref{eq:6}) can be rewritten as 
$$ y(t) = y_0 + \frac{\lambda}{\Gamma(\alpha)} \int_0^t \frac{y(s)}{(t-s)^{1-\alpha}} ds + \frac{1}{\Gamma(\alpha)} \int_0^t \frac{f(s)}{(t-s)^{1-\alpha}} ds. $$
We now apply a Picard-style iteration of the form
$$ y_k(t) = y_0(t)+ \frac{\lambda}{\Gamma(\alpha)} \int_0^t \frac{y_{k-1}(s)}{(t-s)^{1-\alpha}} ds + \frac{1}{\Gamma(\alpha)} \int_0^t \frac{f(s)}{(t-s)^{1-\alpha}} ds, \quad k = 1,2, \cdots $$
where $y_0(t) = y_0, \> \forall t$.

It can be shown that this iteration will converge to (\ref{eq:7}) - see \cite{ref13}. $\square$

\vspace{5mm}

\noindent
\textbf{Lemma 5.}
$$ 1 + \int_0^t \lambda s^{\alpha - 1} E_{\alpha \alpha} (\lambda s^{\alpha}) ds = E_{\alpha} (\lambda s^{\alpha}).$$

\noindent
\textbf{Proof:} Use Definition 1 and integrate the left hand side term by term. $\square$

\vspace{5mm}

\noindent
\textbf{Remark 1.} The function multiplying $f(s)$ in the integrand  of (\ref{eq:7}), namely
$$ G_{\alpha}(t-s) = (t-s)^{\alpha-1} E_{\alpha \alpha} (\lambda (t-s)^{\alpha}),$$
can be viewed as a Green function. For example, when $\alpha = 1$, $G_1(t-s) = e^{\lambda (t-s)}.$

\vspace{5mm}

The generalisation of the class of problems given by (\ref{eq:6}) to the systems case takes the form

\begin{equation}
D_t^{\alpha} y(t) = A y(t) + F(t), \quad y(0) = y_0, \quad y \in \mathbb{R}^m.
\label{eq:9}
\end{equation}
In the case that $F(t) = 0$ the solution of the linear homogeneous system is
\begin{equation}
y(t) = E_\alpha (t^{\alpha} A) y_0.
\label{eq:10}
\end{equation}
From this we can prove

\vspace{5mm}

\noindent
\textbf{Theorem 6.} The solution of (\ref{eq:9}) is given by
\begin{equation}
y(t) = E_{\alpha} (t^{\alpha} A) y_0 + \int_0^t (t-s)^{\alpha-1} E_{\alpha \alpha} ( (t-s)^{\alpha} A) F(s) ds.
\label{eq:11}
\end{equation}

\noindent
\textbf{Proof:} We can use the idea of a Green function. But first of all it is trivial to see from  Definition 1 that
\begin{equation}
\frac{d}{dz} (E_{\alpha} (z^{\alpha} A) ) = A z^{\alpha-1} E_{\alpha \alpha} (z^{\alpha} A).
\label{eq:12}
\end{equation}
Now the solution of (\ref{eq:9}) can be written as
$$y(t) = y_0 + \int_0^t G_{\alpha} (t-s) (Ay_0 + F(s))ds,$$
where $G_{\alpha}(t-s)$ is a matrix Green function, or alternatively
$$
y(t) = \left(I + \int_0^t G_{\alpha}(t-s) A \, ds \right)\,y_0 + \int_0^t G_{\alpha}(t-s)F(s)ds.
$$
Finally it is clear from (\ref{eq:12}) that
$$ G_{\alpha}(t-s) = (t-s)^{\alpha-1} E_{\alpha \alpha} ((t-s)^{\alpha} A) $$
is the Green function and the result is proved by using Lemma 5. $\square$

Note that the proofs of Lemma 4 and Theorem 6 can be found in Podlubny \cite{ref14}.

\section{The solution of mixed index linear systems}

The main focus of this paper is to consider generalisations of (\ref{eq:9}) where there are differing values of $\alpha$ on the left hand side. In its general form, we will let $y^{\top}=(y_1^{\top},\cdots,y_r^{\top}) \in \mathbb{R}^m$ where $y_i \in \mathbb{R}^{m_i}$ and $m = \sum_{i=1}^r m_i$. We will also assume $F(t)^{\top} = (F_1(t)^{\top},\cdots,F_r(t)^{\top})$ and that $A$ can be written in block form $A = (A_{ij})_{i,j=1}^r, \> A_{ij} \in \mathbb{R}^{m_i \times m_j}$. We will also let $\alpha = (\alpha_1,\cdots,\alpha_r)$ and consider a class of linear, non-homogeneous multi-indexed systems of FDEs of the form
\begin{equation}
D_t^{\alpha} y(t) = A y(t) + F(t)
\label{eq:14}
\end{equation}
that we interpret as the system
\begin{equation}
D_t^{\alpha_i}  y_i(t) = \sum_{j=1}^r A_{ij}y_j(t) + F_i(t), \quad i = 1,\cdots,r.
\label{eq:15}
\end{equation}
The index of the system is said to be $r$.

In the case that $F = 0$, then by letting $$E_i = D^{\alpha_i} - A_{1i} $$ we can rewrite (\ref{eq:14}) as
\begin{equation}
M \> y = 0,
\label{eq:1star}
\end{equation}
where $M$ is the block matrix, whose determinant must be zero, with 
\begin{eqnarray*}
M_{ii} &=& E_i, \quad i = 1,\cdots,r \\
M_{ij} &=& -A_{ij}, \quad i \neq j. \\
\end{eqnarray*}

Thus, in the case all $m_{i} = 1$, so that the individual components are scalar and so $m = r$, (\ref{eq:1star}) implies $ \textrm{Det}(M)\> y_r = 0. $

For example, when $r = 2$ this becomes
$$(E_1 E_2 - A_{21}A_{12})\> y_2 = 0$$
or
$$
 (D^{\alpha_1 + \alpha_2} - A_{22}D^{\alpha_1} - A_{11}D^{\alpha_2} + \textrm{Det}(A))\> y_2 = 0.
$$

While, for $r=3$ this gives after some simplification
\begin{eqnarray*}
D^{\alpha_1 + \alpha_2 + \alpha_3} \> y_3 &-& A_{11}D^{\alpha_2 + \alpha_3} \> y_3 - A_{22}D^{\alpha_1+\alpha_3}\> y_3 - A_{33} D^{\alpha_1+\alpha_2} \> y_3   \\
 &+& (A_{22} A_{33} - A_{23}A_{32})D^{\alpha_1}\> y_3 + (A_{11}A_{33} - A_{13}A_{31})D^{\alpha_2}\> y_3  \\
  &+& (A_{11}A_{22} - A_{12}A_{21})D^{\alpha_3}\> y_3 - \textrm{Det}(A) = 0.
\end{eqnarray*}

Clearly there is a general formula for arbitrary $r$ in terms of the cofactors of $A$. In particular, it can be fitted into the framework of  linear sequential FDEs \cite{ref13, ref14, ref15, ref16, ref17}. These take the form
\begin{equation}
D_t^{\beta_0} y_1(t) + \sum_{j=1}^{p} a_j D_t^{\beta_j} y_1(t) = d y_1(t) + f(t), \quad \beta_0 > \beta_1 > \cdots \beta_p.
\label{eq:star}
\end{equation}
However, this characterisation is not particularly simple, useful, or computationally expedient. Furthermore when the $m_i$ are not 1, so that the individual components are not scalar, then there is no simple representation such as (\ref{eq:star}) and new approaches are needed. Before we consider this new approach we note the converse, namely that (\ref{eq:star}) can always be written in the form of (\ref{eq:14}) for a suitable matrix $A$ with a special structure. In particular we can write (\ref{eq:star}) in the form of (\ref{eq:14}) with $p = r-1$ as an $r$ dimensional, $r$ index problem with $\alpha = (\beta_0,\beta_{1},\cdots,\beta_p),$   and
$$ A = \left( \begin{array}{ccccc} d & -a_1 & -a_2 & \cdots  & -a_p \\ 0 & 1 & 0 & \cdots & 0 \\ \vdots & & \ddots  & &  \\ 0 &  & \cdots &  & 1 \end{array} \right), \quad F(t) = (f(t),0,\cdots,0)^{\top}.$$
For completeness we note in the case that $d = 0$ and $f(t) = 0$, an explicit solution to this problem was given in Podlubny \cite{ref14}. This can be found by considering the transfer function (see section 4) given by
$$ H(s) = \frac{1}{s^{\beta_0} + a_1 s^{\beta_1} + \cdots + a_p s^{\beta_p}}. $$
By finding the poles of this function and converting back to the untransformed domain, Podlubny gives the solution as 
\begin{eqnarray*}
y_1(t) &=& \sum_{m=0}^\infty \frac{(-1)^m}{m!} \sum_{\begin{array}{c} k_0+k_1+ \cdots + k_{p-2}=m \\ k_i \geq 0 \\ \end{array} } \left( \begin{array}{c} m \\ k_0 \cdots k_{p-2} \\ \end{array} \right) \prod_{i=0}^{p-2} (a_{p-i})^{k_i} \times \\
&& \epsilon_m (t,-a_1; \beta_0 -\beta_1, \beta_0 + \sum_{j=0}^{p-2} (\beta_1 - \beta_{p-j})k_j + 1) \\
\end{eqnarray*}
where 
\begin{eqnarray*}
\epsilon_k(t,y; \alpha,\beta) &=& t^{k \alpha + \beta - 1} E_{\alpha, \beta}^k (y t^{\alpha}) \\
E_{\alpha, \beta}^k (z) &=& \sum_{i=0}^\infty \frac{(i+k)! \, z^i}{i! \, \Gamma(\alpha(i+k) + \beta)}. \\
\end{eqnarray*}

We now return to the index-2 problem (\ref{eq:1}) and (\ref{eq:sec2eq1}). We first claim that the solution takes the matrix form

\begin{eqnarray}
y_1 & = & \alpha_{00} + \sum_{n=1}^\infty \sum_{j=0}^{n-1} \alpha_{n,j+1} \frac{t^{n \alpha + j (\beta-\alpha)}}{\Gamma(1+n \alpha + j(\beta-\alpha))} z  \nonumber \\
&& \label{eq:17} \\
 y_2 & = & \beta_{00} + \sum_{n=1}^\infty \sum_{j=1}^{n} \beta_{n,j} \frac{t^{n \alpha + j (\beta-\alpha)}}{\Gamma(1+n \alpha + j(\beta-\alpha))} z, \nonumber 
\end{eqnarray}
where the $\alpha_{n,j},\> \beta_{n,j}$ are appropriate matrices, of size $m_1 \times m$ and $m_2 \times m,$ respectively, that are to be determined.

We now use the fact that

\begin{eqnarray}
D_t^{\alpha} \frac{t^{n \alpha + j(\beta-\alpha)}}{\Gamma(1+n \alpha + j(\beta-\alpha))} &=& \frac{1}{\Gamma(1+(n-1)\alpha + j (\beta-\alpha))} \, t^{(n-1)\alpha + j (\beta-\alpha)} \nonumber \\
&& \label{eq:18} \\
D_t^{\beta} \frac{t^{n \alpha + j(\beta-\alpha)}}{\Gamma(1+n \alpha + j(\beta-\alpha))} &=& \frac{1}{\Gamma(1+(n-1)\alpha + (j-1) (\beta-\alpha))} \, t^{(n-1)\alpha + (j-1) (\beta-\alpha)}. \nonumber
\end{eqnarray}
Using (\ref{eq:17}) and (\ref{eq:18}) the left hand side of (\ref{eq:1}) is
\begin{eqnarray*}
D_t^{\alpha} y_1 &=& \sum_{n=1}^\infty \, \sum_{j=0}^{n-1} \alpha_{n,j+1} \frac{t^{(n-1)\alpha + j(\beta-\alpha)}}{\Gamma(1+(n-1)\alpha + j(\beta-\alpha))} \, z \\
D_t^{\beta} y_2 &=& \sum_{n=1}^\infty \, \sum_{j=0}^{n-1} \beta_{n,j+1} \frac{t^{(n-1)\alpha + j(\beta-\alpha)}}{\Gamma(1+(n-1)\alpha + j(\beta-\alpha))} \, z 
\end{eqnarray*}
that can be written in matrix form as 
\begin{equation}
\sum_{n=0}^\infty \, \sum_{j=0}^n \left( \begin{array}{c} \alpha_{n+1,j+1} \\ \beta_{n+1,j+1} \end{array} \right) \frac{t^{n \alpha + j(\beta-\alpha)}}{\Gamma(1+n \alpha + j(\beta-\alpha))} \, z.
\label{eq:19}
\end{equation}
If we define
\begin{equation}
\alpha_{n,n+1} = 0, \quad \beta_{n0} = 0, \quad n=1,2,\cdots
\label{eq:20}
\end{equation}
then the right hand side of (\ref{eq:1}) is
\begin{equation}
A \left( \left( \begin{array}{c} \alpha_{00} \\ \beta_{00} \end{array} \right) + \sum_{n=1}^\infty \, \sum_{j=0}^n \left( \begin{array}{c} \alpha_{n,j+1} \\ \beta_{nj} \end{array} \right) \frac{t^{n \alpha + j(\beta-\alpha)}}{\Gamma(1+n \alpha + j(\beta-\alpha))} \right) \, z.
\label{eq:21}
\end{equation}
Equating (\ref{eq:19}) and (\ref{eq:21}) we find along with (\ref{eq:20}) that for $n=0,1,2,\cdots$
\begin{equation}
\left( \begin{array}{c} \alpha_{00} \\ \beta_{00} \end{array} \right) = I_m, \> \left( \begin{array}{c} \alpha_{n+1,j+1} \\ \beta_{n+1,j+1} \end{array} \right) = A \left( \begin{array}{c} \alpha_{n,j+1} \\ \beta_{nj} \end{array} \right), \quad j=0,1,\cdots,n.
\label{eq:22}
\end{equation}
In order to get a succinct representation of the solution based on (\ref{eq:17}) and (\ref{eq:22}), it will be convenient to write
$$
p_n(t) = \left( \frac{t^{n \alpha}}{\Gamma(1+n \alpha)},\frac{t^{(n-1)\alpha + \beta}}{\Gamma(1+(n-1)\alpha+\beta)}, \cdots,\frac{t^{n\beta}}{\Gamma(1+n\beta)} \right)^{\top} \otimes I_m, \quad n=1,2,\cdots 
$$
so $p_n(t) \in \mathbb{R}^{m(n+1)\times m}$, and let $p_0(t) = I_m$.

We will also define the matrices
\begin{eqnarray*}
L_n &=& \left( \begin{array}{ccccc}
\alpha_{n1} & \alpha_{n2} & \cdots & \alpha_{nn} & 0 \\
0 & \beta_{n1} & \cdots & \beta_{n\, n-1} & \beta_{nn} \\
\end{array} \right) \in \mathbb{R}^{m \times m(n+1)}, \> n=1,2,\cdots \\
L_0 &=& I_m
\end{eqnarray*}
where 0 represents appropriately-sized zero matrices.
Now we note that the recursive relation (\ref{eq:22}) is equivalent to
\begin{equation}
\left( \begin{array}{ccc}
\alpha_{n1}  & \cdots & \alpha_{nn}  \\
 \beta_{n1} & \cdots & \beta_{nn} \\
\end{array} \right) = A\, L_{n-1}, \quad n = 1,2,\cdots .
\label{eq:25}
\end{equation}
Thus we can state the following theorem.

\vspace{5mm}
\noindent
\textbf{Theorem 7.} The solution of the fractional index-2 system $$D_t^{\alpha,\beta}\, y(t) = A\, y(t), \> y(0) = z$$ is given by
\begin{equation}
y(t) = \sum_{n=0}^\infty L_n \, p_n(t) \, z,
\label{eq:26}
\end{equation}
where for $n=1, 2, \cdots$
\begin{eqnarray}
L_n &=& \left( \begin{array}{ccccc}
\alpha_{n1} & \alpha_{n2} & \cdots & \alpha_{nn} & 0 \\
0 & \beta_{n1} & \cdots & \beta_{n\, n-1} & \beta_{nn} \\
\end{array} \right), \quad 
\left( \begin{array}{ccc}
\alpha_{n1} & \cdots & \alpha_{nn}  \\
 \beta_{n1} & \cdots  & \beta_{nn} \\
\end{array} \right) = A \, L_{n-1},  \nonumber \\
L_0 &=& I_m \nonumber \\
p_n(t) &=& \left( \frac{t^{n \alpha}}{\Gamma(1+n \alpha)}, \frac{t^{(n-1)\alpha+\beta}}{\Gamma(1+(n-1)\alpha+\beta)}, \cdots, \frac{t^{n \beta}}{\Gamma(1+n \beta)} \right)^{\top} \otimes I_m. \label{eq:27}
\end{eqnarray}

\vspace{5mm}
\noindent
\textbf{Remarks 2.}

\begin{enumerate}
\item[(i)] In the case $\alpha = \beta$, 
$$p_n(t) = \frac{t^{n \alpha}}{\Gamma(1+n \alpha)} \, (1,\cdots,1)^{\top} \otimes I_m,$$
$$ L_n \, p_n(t) = \frac{t^{n \alpha}}{\Gamma(1+n \alpha)} \sum_{j=1}^n \left( \begin{array}{c} \alpha_{nj} \\ \beta_{nj} \\ \end{array} \right)$$
and with
$$ \sum_{j=1}^n \left( \begin{array}{c} \alpha_{nj} \\ \beta_{nj} \\ \end{array} \right) = A \sum_{j=1}^{n-1} \left( \begin{array}{c} \alpha_{n-1,j} \\ \beta_{n-1,j} \\ \end{array} \right)$$
then (\ref{eq:26}) reduces, as expected, to
$$ y(t) = E_{\alpha} (t^{\alpha} A)z.$$
\item[(ii)] It will be convenient to define the matrix $$P_{\alpha, \beta}(t) = \sum_{n=0}^{\infty} L_n p_n(t)$$ so that the solution (\ref{eq:26}) can be expressed as 
\begin{equation}
y(t) = P_{\alpha, \beta} (t) y_0.
\label{eq:28}
\end{equation}
\item[(iii)] If the fractional index-2 system has initial condition $y(t_0) = z$ then the solution is
\begin{equation}
y(t) = P_{\alpha, \beta} (t-t_0) z.
\label{eq:29}
\end{equation}
\end{enumerate}

We note that in solving (\ref{eq:9}) an equivalent solution to (\ref{eq:11}) is
\begin{eqnarray*}
y(t) &=& E_{\alpha} (t^{\alpha} A) y_0 + I_{\alpha} (G_{\alpha} (t-s) F(s) ) ds, \\
G_{\alpha} (t-s) &=& E_{\alpha} ((t-s)^{\alpha} A), 
\end{eqnarray*}
where $G_{\alpha}$ is the Green function satisfying
\begin{equation}
D_t ^{\alpha} G_{\alpha} (t-s) = A G_{\alpha} (t-s).
\label{eq:31}
\end{equation}
This leads us to give a general result on the solution of the mixed index problem with a time-dependent forcing function
$$D_t^{\alpha, \beta} y(t) = A y(t) + F(t),$$
but first we need the following definition. 

\vspace{5mm}
\noindent
\textbf{Definition 2.}
Let $y(t) = (y_1^{\top} (t), y_2^{\top} (t))^{\top}$, then define
$$I_t^{\alpha, \beta} y(s) ds = \left( I_t^{\alpha} y_1^{\top} (s) ds, I_t^{\beta} y_2^{\top} (s) ds \right) ^{\top}.$$

\vspace{5mm}
\noindent
\textbf{Theorem 8.}
The solution to the fractional index-2 problem
\begin{equation}
D_t^{\alpha, \beta} y(t) = A y(t) + F(t), \quad y(0) = y_0
\label{eq:32}
\end{equation}
is given by
\begin{equation}
y(t) = P_{\alpha, \beta} (t) y_0 + I_t^{\alpha, \beta} \left( P_{\alpha, \beta} (t-s) F(s) ds \right).
\label{eq:33}
\end{equation}

\vspace{3mm}
\noindent
\textbf{Proof:}
The result follows from the above discussion and noting that
$$D_t^{\alpha, \beta} P_{\alpha, \beta} (t) = A P_{\alpha, \beta} (t). \quad \square $$

We now turn to analysing the asymptotic stability of linear fractional index-2 systems.

\section{Asymptotic stability of multi-index systems}

The first contribution to the asymptotic stability analysis of time fractional linear systems was by Matignon \cite{refMatignon}. Given  the linear system $D_t^{\alpha} y(t) = A y(t)$ in Caputo form, then taking the Laplace transform and using the definition of the Caputo derivative gives
$$s^{\alpha} X(s) - s^{\alpha - 1} X(0) = A X(s)$$
or
\begin{equation}
X(s) = \frac{1}{s} (I-s^{-\alpha} A)^{-1} X(0).  \label{eq:3.1}
\end{equation}

\begin{figure}[htb]
\caption{Stability region for single index problem} 
\label{fig:fig1}
\centering
\includegraphics[width=0.65\textwidth]{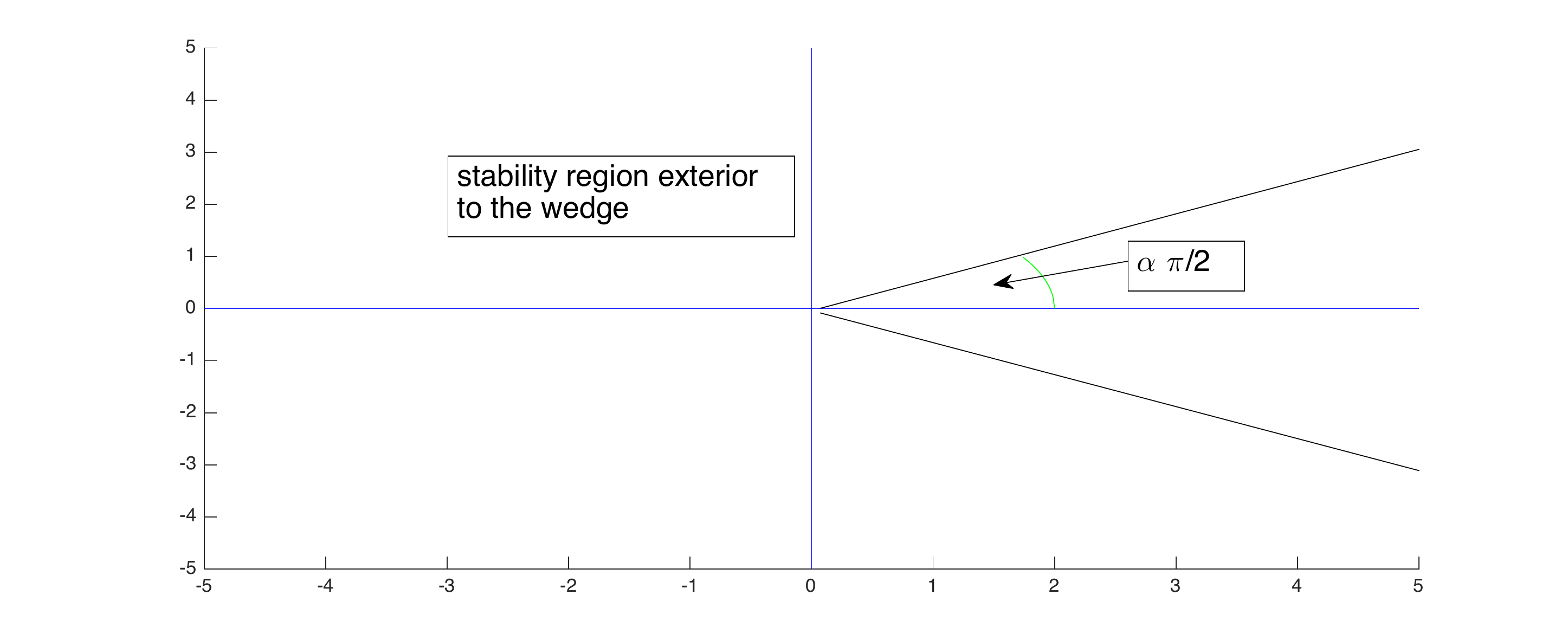} 
\end{figure}

Here $X(s)$ is the Laplace transform of $y(t)$. If we write $w = s^{\alpha},$ then the matrix $s^{\alpha}I - A$ will be nonsingular if $w$ is not an eigenvalue of $A$. In the $w$-domain this will happen if $\mathrm{Re}(\sigma(A)) \leq 0$, where $\sigma(A)$ denotes the spectrum of $A$. In the $s$-domain this will happen if $|\mathrm{Re}(\sigma(A))| \geq \frac{\alpha \pi}{2}.$ That is, the eigenvalues of $A$ lie in the complex plane minus the sector subtended by angle $\alpha \pi$ symmetric about the positive real axis - see Figure \ref{fig:fig1}. 

In fact Laplace transforms are a very powerful technique for studying the asymptotic stability of mixed index fractional systems. Deng et al.  \cite{refDeng} studied the stability of linear time fractional systems with delays using Laplace transforms. Given the delay system
\begin{equation}
\frac{d^{\alpha_i}\,y_i}{dt^{\alpha_i}} = \sum_{j=1}^m a_{ij}y_j(t-\tau_{ij}), \quad i=1,\cdots,m
\label{eq:1-1star}
\end{equation}
then the Laplace transforms results in
\begin{eqnarray}
\Delta(s)\,X &=& b \nonumber \\
\Delta(s) &=& Diag(s^{\alpha_1},\cdots,s^{\alpha_m})  - L   \label{eq:1-2star} \\
L_{ij}& = & a_{ij} e^{-s \tau_{ij}}, \quad i,j = 1,\cdots,m. \nonumber 
\end{eqnarray} 
Hence, Deng et al. proved:

\vspace{5mm}
\noindent
\textbf{Theorem 9.}
If all the zeros of the characteristic polynomial of $\Delta(s)$ have negative real part then the zero solution of (\ref{eq:1-1star}) is asymptotically stable.

\vspace{5mm}

Deng et al. also proved a very nice result in the case that all the indices $\alpha_1,\cdots,\alpha_m$ are rational.

\vspace{5mm}
\noindent
\textbf{Theorem 10.}
Consider (\ref{eq:1-1star}) with no delays and all the $\alpha_i \in (0,1)$ and are rational. In particular let
$$ \alpha_i = \frac{u_i}{v_i}, \quad \textrm{gcd}(u_i,v_i) = 1$$
and let $M$ be the lowest common multiple of all the denominators and set $\gamma = \frac{1}{M}$. Then the problem will be asymptotically stable if all the roots, $\lambda$, of 
$$ p(\lambda) = \textrm{Det}(D-A) = 0, \quad D = \textrm{diag}(\lambda^{Ma_1},\cdots,\lambda^{Ma_m})$$
satisfy $|\arg(\lambda)| > \gamma \frac{\pi}{2}.$ 

\vspace{5mm}
\noindent
\textbf{Remarks 3.}
\begin{enumerate}[label=(\roman*)]
\item If $\alpha_i = \alpha, \> i=1,\cdots,m$ then Theorem 10 reduces to the result of Matignon. The proof of Theorem 10 comes immediately from (\ref{eq:1-2star}) where $p(\lambda)$ is the characteristic polynomial of $\Delta(s)$.
\item A nice survey on the stability (both linear and nonlinear) of fractional differential equations is given in Li and Zhang \cite{refLi}, while Saberi Najafi et al.  \cite{refSaberi} has extended some of these stability results to distributed order fractional differential equations with respect to an order density function. Zhang et al  \cite{refZhang} consider the stability of nonlinear fractional differential equations.
\item Radwan et al.  \cite{refRadwan} note that the stability analysis of mixed index problems reduces to the study of the roots of the characteristic equation 
\begin{equation}
\sum_{i=1}^m \theta_i\,s^{\alpha_i} = 0, \quad 0 < \alpha_i \leq 1.
\label{eq:1-3star}
\end{equation}
\end{enumerate}

In the case that the $\alpha_i$ are arbitrary real numbers, the study of the roots of (\ref{eq:1-3star}) is difficult. By letting $s = e^z$, we can cast this in the framework of quasi (or exponential) polynomials (Rivero et al. \cite{refRivero}). The zeros of exponential polynomials have been studied by Ritt  \cite{Ritt1929}. 

The general form of an exponential polynomial with constant coefficients is
$$f(z) = \sum_{j=0}^k a_j e^{\alpha_j z}. $$
An analogue of the fact that a polynomial of degree $k$ can have up to $k$ roots is expressed by a Theorem due to Tamarkin, P\'{o}lya and Schwengler (see \cite{Ritt1929}).

\vspace{5mm}
\noindent
\textbf{Theorem 11.} Let $P$ be the smallest convex polygon containing the values $\alpha_1,\cdots,\alpha_k$ and let the sides of  $P$ be $s_1,\cdots,s_k$. Then there exist $k$ half strips with half rays parallel to the outer normal to $b_i$ that contain all the zeros of $f$. If $|b_i |$ is the length of $b_i$, then the number of zeros in the $i^{th}$ half strip with modulus less than or equal to $r$ is asymptotically $\frac{r\,|b_i|}{2\pi}.$

\vspace{5mm}

If the $\alpha_i$ are rational and with $M$ the lowest common multiple of the denominators, this reduces to the polynomial
$$\sum_{i=1}^M \theta_i W^i = 0, \quad W = s^{\frac{1}{M}}.$$

This  leads us to think about stability from a control theory point of view. Thus given the system
\begin{equation}
\sum_{j=0}^n a_j D^{\alpha_j} y = \sum_{j=0}^M b_j D^{\beta_j} y
\label{eq:1-4star}
\end{equation}
where
$$ \alpha_n > \cdots > \alpha_0, \quad \beta_M > \cdots \beta_0$$
then the solution of (\ref{eq:1-4star}) can be written in terms of the transfer function 
\begin{equation}
G(s) = \frac{ \sum_{j=0}^m b_j s^{\beta_j}}{\sum_{j=0}^n a_j s^{\alpha_j}} := \frac{Q(s)}{P(s)},
\label{eq:1-5star}
\end{equation}
where $s$ is the Laplace variable (see Rivero et al.  \cite{refRivero}, Petras  \cite{refPetras}).

In the case of the so-called commensurate form in which
$$ \alpha_k = k \alpha, \quad \beta_k = k \beta,$$
then
\begin{equation}
G(s) = \frac{ \sum_{k=0}^m b_k(s^\beta)^k}{\sum_{k=0}^n a_k(s^\alpha)^k} := \frac{Q(s^\beta)}{P(s^\alpha)}.
\label{eq:1-6star}
\end{equation}

Clearly, if $\frac{\beta}{\alpha}$ is rational with $\alpha \geq \beta$ and
$$ \beta = \frac{q}{p} \alpha, \quad q, p \in \mathbb{Z}^+, \quad w = s^{\frac{\alpha}{p}}$$
then (\ref{eq:1-6star}) can be written as
$$G(w) := \frac{Q(w^q)}{P(w^p)}, \quad p,q \in \mathbb{Z}^+, \quad q \leq p.$$

C\v{e}rm\'{a}k and Kisela  \cite{refCermak} considered the specific problem
\begin{equation}
D^\alpha y + a D^\beta y + by = 0, \quad y \in  \mathbb{R},
\label{eq:1-7star}
\end{equation}
where $\alpha = pK, \> \beta = q K, K$ real $\in (0,1), p, q \in \mathbb{Z}^+, p \geq q.$ In this case the appropriate stability polynomial is $P(\lambda) := \lambda^p + a \lambda^q + b$, where $\lambda = s^K$. Based on Theorem 10, (\ref{eq:1-7star}) is asymptotically stable if all the roots of $P(\lambda)$ satisfy $|\arg(\lambda)| > K \frac{\pi}{2}$.

By setting $\lambda = r e^{iK\frac{\pi}{2}}$ and substituting into $P(\lambda) = 0$ and equating real and imaginary parts, it is easily seen that 
\begin{eqnarray*}
r^p \cos \frac{pK\pi}{2} + a\, r^q \cos \frac{qK\pi}{2} + b &=& 0 \\
r^p \sin \frac{pK\pi}{2} + a \, r^q \sin \frac{qK\pi}{2} &=& 0.
\end{eqnarray*}

This leads to the following result, given in C\v{e}rm\'{a}k and Kisela \cite{refCermak}.

\vspace{5mm}
\noindent
\textbf{Theorem 12.}
Equation (\ref{eq:1-7star}) is asymptotically stable with $\alpha > \beta > 0$ real and $\frac{\alpha}{\beta}$ rational if
\begin{eqnarray*}
\beta &<& 2, \quad \alpha - \beta < 2 \\
b &>& 0, \quad a > \frac{-\sin \frac{\alpha \pi}{2}}{(\sin \frac{\beta \pi}{2})^{\frac{\beta}{\alpha}}\, (\sin \frac{(\alpha - \beta) \pi}{2})^{\frac{\alpha-\beta}{\alpha}}} \, b^{\frac{\alpha-\beta}{\alpha}}. \\
\end{eqnarray*}

We now follow this idea but for arbitrarily sized systems in our mixed index format, and this leads to slight modifications to (\ref{eq:1-7star}). We first make a slight simplification and take $m_1 = m_2$ and we also assume that $A_2$ is nonsingular, then problem (\ref{eq:1}) leads to 
$$y_2 = A_2^{-1}\,(D^\alpha\,I - A_1)\, y_1$$
and substituting into the equation for $y_1$ gives  
$$ \begin{array}{c}
(D^{\alpha + \beta}\,I - B_2 D^\alpha \, I - \bar{A}_1 D^\beta \, I + B_2 \bar{A}_1 - B_1 A_2)\,A_2^{-1}\,y_1 = 0 \\
\bar{A}_1 = A_2^{-1}A_1 A_2. \\
\end{array}$$
This leads us to consider the roots of the characteristic function
\begin{equation}
P(\lambda) := \textrm{Det}(D^{\alpha+\beta}I - B_2 D^\alpha I - \bar{A}_1 D^\beta I + B_2 \bar{A}_1 - B_1 A_2) = 0.
\label{eq:1-8star}
\end{equation}
In the scalar case this gives an extension to (\ref{eq:1-7star}) where the characteristic equation is
\begin{equation}
P(\lambda) = \lambda^{\alpha+\beta} - B_2 \lambda^\alpha - A_1 \lambda^\beta + \textrm{Det}(A).
\label{eq:1-9star}
\end{equation}

Now reverting to Laplace transforms of (\ref{eq:1}) and (\ref{eq:sec2eq1}) then 
\begin{eqnarray*}
s^{\alpha} X_1(s) - s^{\alpha - 1}X_1(0) &=& A_1 X_1(s) + A_2 X_2(s) \\
s^{\beta} X_2(s) - s^{\beta - 1}X_2(0) &=& B_1 X_1(s) + B_2 X_2(s). 
\end{eqnarray*}

This can be written in systems form as
\begin{equation} (D_1 - A) X(s) = D_2 X(0),
\label{eq:3.2}
\end{equation}
where
$$ D_1 = \left( \begin{array}{cc} s^{\alpha} I & 0 \\ 0 & s^{\beta}I \end{array} \right), \quad
D_2 = \left( \begin{array}{cc} s^{\alpha-1} I & 0 \\ 0 & s^{\beta-1}I \end{array} \right) $$
or alternatively as
\begin{equation}
X(s) = \frac{1}{s} (I-D_1^{-1}A)^{-1}X(0).
\label{eq:3.2b}
\end{equation}

This can now be considered as a generalised eigenvalue problem. From (\ref{eq:3.2}) we require $D_1 - A$ to be nonsingular. That is
$$  \left( \begin{array}{cc} s^{\alpha} I-A_1 & -A_2 \\ -B_1 & s^{\beta}I-B_2 \end{array} \right)\, v = 0 \implies v = 0.$$

Let us write $v = (v_1^{\top},\, v_2^{\top})^{\top}$ and assume $\alpha \geq \beta$ and that $s^{\beta}I - B_2$ is nonsingular, so that from the previous analysis this means
\begin{equation}
|\mathrm{Re}(\sigma(B_2))| \geq \frac{\beta \pi}{2}.
\label{eq:3.4}
\end{equation}
Hence 
\begin{eqnarray*}
v_2 &=& (s^{\beta}I - B_2)^{-1} B_1 v_1 \\
( (s^{\alpha}I - A_1) - A_2(s^{\beta}I - B_2)^{-1} B_1) v_1 &=& 0.
\end{eqnarray*}
Thus  (\ref{eq:3.4}) and
\begin{equation}
\textrm{Det}((s^{\alpha}I - A_1) - A_2(s^{\beta}I - B_2)^{-1}B_1) = 0
\label{eq:3.5}
\end{equation}
define the asymptotic stability boundary - see also (\ref{eq:1-8star}).

In order to make this more specific, let $m_1 = m_2 = 1$ and 
\begin{equation}
A = \left[ \begin{array}{cc} d & b \\ a & d \end{array} \right], \quad d < 0.
\label{eq:eq44a}
\end{equation}

Note that $\sigma(A) = \{d \pm \sqrt{ab}\}.$ Then (\ref{eq:3.5}) becomes
\begin{equation}
(s^{\alpha}-d)(s^{\beta}-d) - ab = 0.
\label{eq:3.6}
\end{equation}
Furthermore, let $b = -a = \theta$, so that the eigenvalues of $A$ are $d \pm i \theta$ and (\ref{eq:3.6}) becomes
\begin{equation}
(s^{\alpha} - d)(s^{\beta} - d) + \theta^2 = 0.
\label{eq:3.7}
\end{equation}
If we now assume that $$s = r e^{i \frac{\pi}{2}},$$
which defines the asymptotic stability boundary (the imaginary axis) when $\alpha = \beta = 1,$ then (\ref{eq:3.7}) becomes
\begin{equation}
\theta^2 = -(r^{\alpha} e^{i \frac{\pi \alpha}{2}} - d)\,(r^{\beta} e^{i \frac{\pi \beta}{2}} - d).
\label{eq:3.8}
\end{equation}
Now since $\theta$ and $d$ are real, the imaginary part of the right hand side of (\ref{eq:3.8}) must be zero, so that
\begin{equation} 
r^{\alpha+\beta} \sin \frac{\alpha+\beta}{2} \pi = d (r^{\alpha} \sin \frac{\alpha \pi}{2} + r^{\beta} \sin \frac{\beta \pi}{2} ).
\label{eq:3.9}
\end{equation}

Hence
\begin{equation}
-\theta^2 = r^{\alpha+\beta} \cos \frac{\alpha+\beta}{2} \pi - d (r^{\alpha} \cos \frac{\alpha \pi}{2} + r^{\beta} \cos \frac{\beta \pi}{2}) + d^2.
\label{eq:3.10}
\end{equation}

Equations (\ref{eq:3.9}) and (\ref{eq:3.10}) will define the asymptotic stability boundary with $\theta$ as a function of $d$. Rewriting (\ref{eq:3.9}) as
\begin{equation}
d = \frac{r^{\alpha+\beta} \sin \frac{\alpha+\beta}{2} \pi}{r^\alpha \sin \frac{\alpha \pi}{2} + r^\beta \sin \frac{\beta \pi}{2} }.
\label{eq:star2}
\end{equation}
and substituting (\ref{eq:3.10}) leads after simplification to
$$ \frac{\theta^2}{d^2} = \frac{1}{r^{\alpha+\beta}(\sin \frac{\alpha+\beta}{2} \pi)^2} 
\left[ \sin \frac{\alpha+\beta}{2} \pi ( \frac{r^{2 \alpha}}{2} \sin \alpha \pi + \frac{r^{2 \beta}}{2} \sin \beta \pi) \right. $$
$$ \left. \quad \quad \quad  \quad \quad \quad - \cos \frac{\alpha+\beta}{2} \pi (r^{2 \alpha} \sin^2 \frac{\alpha \pi}{2} + r^{2 \beta} \sin^2 \frac{\beta \pi}{2} + 2 r^{\alpha+\beta} \sin \frac{\alpha \pi}{2} \sin \frac{\beta \pi}{2}) \right].
$$
Using the relationships
\begin{eqnarray*}
\sin^2 \theta &=& \frac{1}{2}(1-\cos 2 \theta) \\
\sin A \sin B + \cos A \cos B &=& \cos (A-B)
\end{eqnarray*}
gives
\begin{eqnarray}
\frac{\theta^2}{d^2} &=& \frac{1}{2r^{\alpha+\beta} \sin^2 \frac{\alpha+\beta}{2} \pi} \left( (r^{2 \alpha} + r^{2 \beta}) \, ( \cos \frac{\alpha-\beta}{2} \pi - \cos \frac{\alpha+\beta}{2} \pi) \right. \nonumber \\
&& \left. \quad \quad \quad \quad \quad \quad  \quad \quad \quad  - 4r^{\alpha+\beta} \sin \frac{\alpha \pi}{2} \sin \frac{\beta \pi}{2} \cos \frac{\alpha+\beta}{2} \pi \right).
\label{eq:star3}
\end{eqnarray}
Since 
$$ \cos \frac{\alpha - \beta}{2} \pi - \cos \frac{\alpha+\beta}{2} \pi = 2 \sin \frac{\alpha \pi}{2} \sin \frac{\beta \pi}{2}$$
and letting $x = r^{\alpha - \beta}$, then we can write (\ref{eq:star3}) as
\begin{equation}
\left( \frac{\theta}{d}\right)^2 = \frac{ \sin \frac{\alpha \pi}{2} \sin \frac{\beta \pi}{2}}{\sin^2 \frac{\alpha+\beta}{2} \pi} \left( \frac{x^2 + 1}{x} - 2 \cos \frac{\alpha + \beta}{2} \pi \right).
\label{eq:2-1star}
\end{equation}
Furthermore, we can write (\ref{eq:star2}) as
\begin{equation}
d = \frac{x^{\frac{\alpha}{\alpha-\beta}} \sin \frac{\alpha+\beta}{2} \pi}{x \sin \frac{\alpha \pi}{2} + \sin \frac{\beta \pi}{2}}.
\label{eq:2-2star}
\end{equation}
It is easily seen that as a function of $x$ the minimum of (\ref{eq:2-1star}) is when $x = 1$. Thus
\begin{eqnarray*}
\frac{\theta}{d} & \geq & \frac{\sqrt{2 \sin \frac{\alpha \pi}{2} \sin \frac{\beta \pi}{2}}}{\sin \frac{\alpha+\beta}{2} \pi} \, \sqrt{1 - \cos \frac{\alpha+\beta}{2} \pi} \\
&=& \frac{2 \sqrt{\sin \frac{\alpha \pi}{2} \sin \frac{\beta \pi}{2}} \, \sin \frac{\alpha+\beta}{4} \pi}{2 \sin \frac{\alpha+\beta}{4} \pi \cos \frac{\alpha+\beta}{4} \pi} \\
&=& \frac{\sqrt{\sin \frac{\alpha \pi}{2} \sin \frac{\beta \pi}{2}}}{\cos \frac{\alpha+\beta}{4} \pi}.
\end{eqnarray*}
Thus we have proved the following result.

\vspace{5mm}
\noindent
\textbf{Theorem 13.}
Given the mixed index problem with $A$ as in (\ref{eq:eq44a}), the angle for asymptotic stability $\hat{\theta} = \arctan(\frac{\theta}{d})$ satisfies
\begin{equation}
\tan \hat{\theta} \in \left[ \frac{\sqrt{\sin \frac{\alpha \pi}{2} \sin \frac{\beta \pi}{2}}}{\cos \frac{\alpha+\beta}{4} \pi}, \infty \right),
\label{eq:2-3star}
\end{equation}
or in radians with $\tilde{\theta} = \frac{1}{\pi} \arctan(\frac{\theta}{d})$
$$
\tilde{\theta} \in \frac{1}{\pi} \left[ \arctan \frac{\sqrt{\sin \frac{\alpha \pi}{2} \sin \frac{\beta \pi}{2}}}{\cos \frac{\alpha+\beta}{4} \pi}, \arctan \frac{\pi}{2} \right]
$$
with the minimum occuring with
\begin{equation}
d = \frac{\sin \frac{\alpha+\beta}{2} \pi}{\sin \frac{\alpha \pi}{2} + \sin \frac{\beta \pi}{2}}.
\label{eq:2-4star}
\end{equation}

\vspace{5mm}
\noindent
\textbf{Remarks 4.}
We have the following results for $\hat{\theta}$ in three particular cases:
\begin{enumerate}[label=(\roman*)]
\item $\alpha = \beta: \quad \hat{\theta} = \alpha \frac{\pi}{2}$, since in this case $(\frac{\theta}{d})^2 = \tan^2 \frac{\alpha \pi}{2}$.
\item $\alpha + \beta = 1: \quad \hat{\theta} \in (\sqrt{\sin \alpha \pi}, \frac{\pi}{2}), \> \alpha \in [ \frac{1}{2}, 1]$. In the case $\alpha + \beta = 1$ we see from (\ref{eq:2-1star}) that 
$$ \left( \frac{\theta}{d}\right)^2 = \sin \alpha \pi \left( \frac{x^2 + 1}{2x} \right).$$
Letting $\alpha = \frac{1}{2} + \epsilon$ with $\epsilon > 0$ small, then $x = r^{2 \epsilon}$. This means that $\frac{x^2+1}{2x}$, as a function of $r$, is very shallow apart from when $r$ is near the origin or very large. Hence the asymptotic stability boundary will be almost constant over long periods of $d$ when $\alpha$ and $\beta$ are close together.
\item $\alpha = 2 \beta: \quad \hat{\theta} \in [ \frac{ \sin \frac{\beta \pi}{2} \sqrt{2 \cos \frac{\beta \pi}{2}}}{\cos \frac{3 \beta \pi}{4}}, \frac{\pi}{2}), \> \beta \in (0,\frac{1}{2}].$
\end{enumerate}

Letting 
$$K = \frac{\sin \frac{\alpha \pi}{2} \sin \frac{\beta \pi}{2}}{\sin^2 \frac{\alpha + \beta}{2} \pi}, \quad
L = 2 \cos \frac{\alpha+\beta}{2} \pi, \> \phi = \frac{\theta}{d},$$
 we can write (\ref{eq:2-1star}) and (\ref{eq:2-2star}) as
\begin{equation}
x^2 - x(L+\frac{\phi^2}{K})x + 1 = 0
\label{eq:2-5star}
\end{equation}
\begin{equation}
x^{\frac{\alpha}{\alpha - \beta}} - x\,d_\alpha - d_\beta = 0,
\label{eq:2-6star}
\end{equation}
where
$$d_\alpha = d \frac{\sin \frac{\alpha \pi}{2}}{\sin \frac{\alpha+\beta}{2} \pi}, \quad d_\beta = d \frac{\sin \frac{\beta \pi}{2}}{\sin \frac{\alpha+\beta}{2} \pi}. $$
Due to the nonlinearities in (\ref{eq:2-6star}) it is hard to determine an explicit simple relation between $\phi$ and $d$ except if $\alpha = 2 \beta$. In this case we make use of the following Lemma.

\vspace{5mm}
\noindent
\textbf{Lemma 14.}
If $x^2-ax+b = 0$ and $x^2-cx+d = 0$ then there is a solution 
\begin{equation}
\begin{array}{ll}
x = 0, & b = d  \\
x^2 - ax + b = 0, & a = c, b = d  \\
x = \frac{d-b}{c-a}, & c \neq a \textrm{   and   } (d-b)^2 = (c-a)(ad-bc).
\end{array}
\label{eq:2-7star}
\end{equation}

\vspace{2mm}
\noindent
\textbf{Proof:}
By subtraction of the two equations and substitution.
$\square$

\vspace{5mm}
In the case of (\ref{eq:2-5star}) and (\ref{eq:2-6star}) then (\ref{eq:2-7star}) becomes
$$(1+d_\beta)^2 = (P-d_\alpha)(Pd_\beta + d_\alpha), \quad P = L + \frac{\phi^2}{K},$$
that is
$$P^2 d_\beta - P d_\alpha(d_\beta -1) - (d_\alpha^2 + (1+d_\beta)^2) = 0.$$
Hence
\begin{equation}
2 d_\beta P = d_\alpha(d_\beta - 1) \pm (1 + d_\beta) \sqrt{d_\alpha^2 + 4 d_\beta}.
\label{eq:2-8star}
\end{equation}
Note that
$$\phi^2 = K\,P - K\, L$$
and $$d_\alpha d_\beta = d^2\,K.$$
Some manipulation from (\ref{eq:2-8star}) leads to
$$\phi^2 = \frac{1}{2} \left(\frac{d_\alpha}{d}\right)^2\,\left(d_\beta - 1 \pm (1+d_\beta) \sqrt{1+4 \frac{d_\beta}{d^2_\alpha}} - 2L \frac{d_\beta}{d_\alpha}\right).$$
Now since $\alpha = 2 \beta$, this reduces to
\begin{eqnarray}
\phi^2 &=& \frac{1}{2} \left( \frac{\sin \beta \pi}{\sin \frac{3 \beta}{2} \pi}\right)^2 \left(d_\beta - 1 \pm (1+d_\beta) \sqrt{1 + \frac{4}{d} \frac{\sin \frac{\beta}{2} \pi \sin \frac{3 \beta}{2} \pi}{(\sin \beta \pi)^2}} - 2 \frac{\cos \frac{3\beta}{2} \pi}{\cos \frac{\beta}{2} \pi} \right) \nonumber \\
d_\beta  &=& d \frac{\sin \frac{\beta}{2} \pi}{\sin \frac{3\beta}{2} \pi}. 
\label{eq:2-9star}
\end{eqnarray}
By taking $\tilde{\theta} = \arctan(\phi)$ this gives an explicit relationship between $\tilde{\theta}$ and $d$ for the case $\alpha = 2\beta$.

\vspace{5mm}
\noindent
\textbf{Remarks 5.}
\begin{itemize}
\item $\beta = \frac{1}{2}, \> \alpha = 1$ gives 
\begin{equation}
\tan \tilde{\theta} = \sqrt{(1+d)(1+\sqrt{1+\frac{2}{d}})}.
\label{eq:2-10star}
\end{equation}
\item $\beta = \frac{1}{3}, \> \alpha = \frac{2}{3}$ gives 
\begin{equation}
\tan \tilde{\theta} = \sqrt{\frac{3}{8}} \sqrt{(1+\frac{d}{2}) \sqrt{1+\frac{8}{3d}} + \frac{d}{2} - 1}.
\label{eq:2-11star}
\end{equation}
\end{itemize}

It is clear from (\ref{eq:2-9star}) that when $d=0$ and $d=\infty$, then $\theta = \frac{\pi}{2}$ and then the angle will make an excursion from $\frac{\pi}{2}$ down to a minimum value and back to $\frac{\pi}{2}$ as $d$ increases. For example, in the case of $\beta = \frac{1}{2}, \> \alpha = 1$ we can show from (\ref{eq:2-10star}) that the minimum value of the angle is when
$$d = \sqrt{2} - 1, \quad \tan \tilde{\theta} = \sqrt{ \sqrt{2} + \sqrt{4 + 3 \sqrt{2}}}.$$
Some of these aspects are shown in the Simulations and Results section.

\section{Further analysis}

Returning to (\ref{eq:3.2}) and taking $m_1 = m_2 = 1$ and 
$$ A = \left( \begin{array}{rr}
a_1 & a_2 \\ b_1 & b_2 \\
\end{array} \right) $$
then the Laplace transform in (\ref{eq:3.2b}) is
\begin{equation}
X(s) = \frac{1}{\textrm{Det}(s)} \left( s^{\alpha+\beta-1} X(0) + \left( \begin{array}{rr} a_2 \\ -a_1 \\ \end{array} \right) s^{\beta-1} X_2(0) + \left( \begin{array}{rr} -b_2 \\ b_1 \\ \end{array} \right) s^{\alpha-1} X_1(0) \right)
\label{eq:eq65}
\end{equation}
where
\begin{eqnarray*}
\textrm{Det}(s) &=&  s^{\alpha+\beta} - a_1 s^{\beta} -b_2 s^{\alpha} + D_A,  \\
D_A &=& a_1 b_2 - a_2 b_1 = \textrm{Det}(A).
\end{eqnarray*}

Now if $\alpha$ and $\beta$ are rational $(\alpha \leq \beta)$ 
$$ \alpha = \frac{m}{n}, \quad \beta = \frac{p}{q}, \quad m \leq n, \> p \leq q, \quad \textrm{positive integers}$$
and with $z = s^{\frac{1}{nq}},$ then
\begin{equation}
\textrm{Det}(z) =z^{mq+np} - a_1 z^{np} - b_2 z^{mq} + D_A.
\label{eq:paper2eq6}
\end{equation}
Hence (\ref{eq:eq65}) gives
\begin{equation}
X_1(z) = \frac{1}{z^{(n-m)q}\textrm{Det}(z)} \left( (z^{np}-b_2)X_1(0) + a_2 z^{np-mq}X_2(0) \right)
\label{eq:paper2eq7}
\end{equation}
\begin{equation}
X_2(z) = \frac{1}{z^{(n-m)q}\textrm{Det}(z)} \left( b_1X_1(0) + (z^{np}-a_1z^{np-mq})X_2(0) \right).
\label{eq:paper2eq8}
\end{equation}
From Descartes rule of sign, then (\ref{eq:paper2eq6}) will have at most 4 real zeros if $mq+np$ is even, and at most 5 real zeros if $mq+np$ is odd. 

Now factorise
$$
\textrm{Det}(z) = \Pi_{j=1}^N (z-\lambda_j), \quad  N = mq + np,
$$
where there are at most 4 real zeros if $N$ is even and at most 5 real zeros if $N$ is odd. Then 
using (\ref{eq:paper2eq7}) and (\ref{eq:paper2eq8}) we can write
$$ X_i(s) = \frac{s^{\frac{1}{nq}}}{s^{1-\alpha+\frac{1}{nq}}} \, \sum_{j=1}^N \frac{A_j^{(i)}}{s^{\frac{1}{nq}}-\lambda_j}, \quad i=1, \> 2$$
where the $A_j^{(i)}$ can be found by writing
$$ \frac{p_i(z)}{\textrm{Det}(z)} = \sum_{j=1}^N \frac{A_j^{(i)}}{z-\lambda_j}, \quad i=1,2$$
where 
\begin{eqnarray*}
p_1(z) &=& X_1(0) z^{np} + X_2(0) a_2 z^{np-mq} - b_2 X_1(0) \\
p_2(z) &=& X_2(0) z^{np} - X_2(0) a_1 z^{np-mq} + b_1 X_1(0).
\end{eqnarray*}

Using Lemma 2 with 
$$ \tilde{\alpha} = \frac{1}{nq}, \quad \tilde{\beta} = 1 - \alpha + \tilde{\alpha}$$
leads to the following result.

\vspace{5mm}
\noindent
\textbf{Theorem 15.} The solution of the mixed index 2 problem  with $\alpha = \frac{m}{n}$, $\beta = \frac{p}{q}, \> m \leq n, \> p \leq q$ all positive integers is, with $N = mq+np$, given by
\begin{eqnarray}
y(t) &=& \sum_{j=1}^N A_j E_{\frac{1}{nq},1-\alpha+\frac{1}{nq}} (\lambda_j t^{\frac{1}{nq}})
\label{eq:paper2eq14} \\
A_j &=& (A_j^{(1)},A_j^{(2)})^\top, \nonumber
\end{eqnarray}
where the $\lambda_j$ are the zeros of (\ref{eq:paper2eq6}) and the $A_j$ are the coefficients in the partial fraction expansion.

\vspace{2mm}
\noindent
\textbf{Remarks 6.} 
\begin{enumerate}[label=(\roman*)]
\item In the case that $\alpha = \beta$ then (\ref{eq:paper2eq14}) should collapse to the solution
\begin{equation}
y(t) = E_{\alpha} (t^\alpha A) y(0),
\label{eq:paper2eq15}
\end{equation}
and this is not immediately clear. However, in this case, $mq = np$ and so
$$D(z) = z^{2np} - (a_1 + b_2)z^{np} + D(A)$$
which is a quadratic function in $z^{np}$ while the equivalent $p_1$ and $p_2$ numerator functions are linear in $z^{np}$. Thus in (\ref{eq:paper2eq14}) $N$ is replaced by 2,  $\frac{1}{nq}$ is replaced by $\alpha$, and $1-\alpha+\frac{1}{nq}$ becomes 1. Thus (\ref{eq:paper2eq14}) reduces to
$$ y(t) = \sum_{j=1}^2 A_j E_{\alpha}(\lambda_j t^\alpha)$$
that then becomes (\ref{eq:paper2eq15}).
\item In the case that $\alpha$ is rational and $\beta = K \alpha, \> K \> \textrm{a positive integer},$ then
\begin{equation}
\textrm{Det}(s) = (s^\alpha)^{K+1} - a_1 (s^\alpha)^K - b_2 s^\alpha + D_A.
\label{eq:paper2eq16}
\end{equation}
If we factorise
$$\textrm{Det}(s) = \Pi_{j=1}^{K+1} (s^\alpha - \lambda_j)$$
and find $A_j^{(1)}, \> A_j^{(2)}, \> j=1,\cdots,K+1$ by
\begin{equation}
\sum_{j=1}^{K+1} A_j \frac{1}{s^\alpha - \lambda_j} = \frac{1}{\textrm{Det}(s)} \left(  (s^\alpha)^K X(0) + \left( \begin{array}{r} a_2 \\ -a_1 \end{array} \right) (s^\alpha)^{K-1} X_2(0) + \left( \begin{array}{r} -b_2 \\ b_1 \end{array} \right) X_1(0) \right)
\label{eq:paper2eq17}
\end{equation}
then we have the following Corollary.
\end{enumerate}

\vspace{5mm}
\noindent
\textbf{Corollary 16.} The solution of the mixed index 2 problem with $\alpha$ rational and  $\beta  = K \alpha$, $K$ a positive integer, is given by
$$ y(t) = \sum_{j=1}^{K+1} A_j E_{\alpha} (\lambda_j t^\alpha),$$
where the vectors $A_j$ and ``eigenvalues" $\lambda_j$ satisfy (\ref{eq:paper2eq17}).

As a particular example, take $K = 2, \> \alpha = \frac{p}{q},$ then the $\lambda_j$ and $A_j$ in Corollary 16 satisfy
$$ D(z): = \Pi_{j=1}^3 (z-\lambda_j) := z^3 - a_1 z^2 - b_2 z + D_A = 0$$
and
$$\sum_{j=1}^3 A_j \frac{1}{z-\lambda_j} = \frac{1}{D(z)} \left( X_0\, z^2 + \left( \begin{array}{r} a_2 \\ -a_1 \end{array} \right) X_2(0) z + \left( \begin{array}{r} -b_2 \\ b_1 \end{array} \right) X_1(0) \right).$$
In other words
\begin{eqnarray*}
A_1 (z-\lambda_2)(z-\lambda_3) + A_2 (z-\lambda_1)(z-\lambda_3) + A_3 (z-\lambda_1)(z-\lambda_2) && \\
= X_0 z^2 + \left( \begin{array}{r} a_2 \\ -a_1 \end{array} \right) X_2(0) z + \left( \begin{array}{r} -b_2 \\ b_1 \end{array} \right) X_1(0) &&
\end{eqnarray*}
or
$$[A_1 \>\> A_2 \>\> A_3] = \left[X_0,\> \left( \begin{array}{r} a_2 \\ -a_1 \end{array} \right) X_2(0), \> \left( \begin{array}{r} -b_2 \\ b_1 \end{array} \right) X_1(0) \right] \> S^{-1}$$
with
$$ S = \left[ \begin{array}{rrr} 1 & -(\lambda_2+\lambda_3) & \lambda_2 \lambda_3 \\
1 & -(\lambda_1+\lambda_3) & \lambda_1 \lambda_3 \\ 1 & -(\lambda_1+\lambda_2) & \lambda_1 \lambda_2 \end{array} \right]. $$

Clearly in the case described by Corollary 16, writing the solution as a linear combination of generalised Mittag-Leffler functions makes the evaluation of the solution much more computationally efficient.

\section{Simulations and results}


\begin{figure}[tbhp]
\captionsetup{width=0.75\textwidth}
\caption{Stability region, above the blue line, for choosing $d$ and $\theta$, when the eigenvalues of $A$ are $d \pm i \theta$, $\alpha = \frac{1}{2}, \beta = 1.$ The logarithmic scale is explored in the right hand figure where the stability boundary dips below the angle $\frac{3\pi}{8}$.} 
\label{fig:fig2}
\centering
$$\begin{array}{cc}
\includegraphics[width=0.35\textwidth]{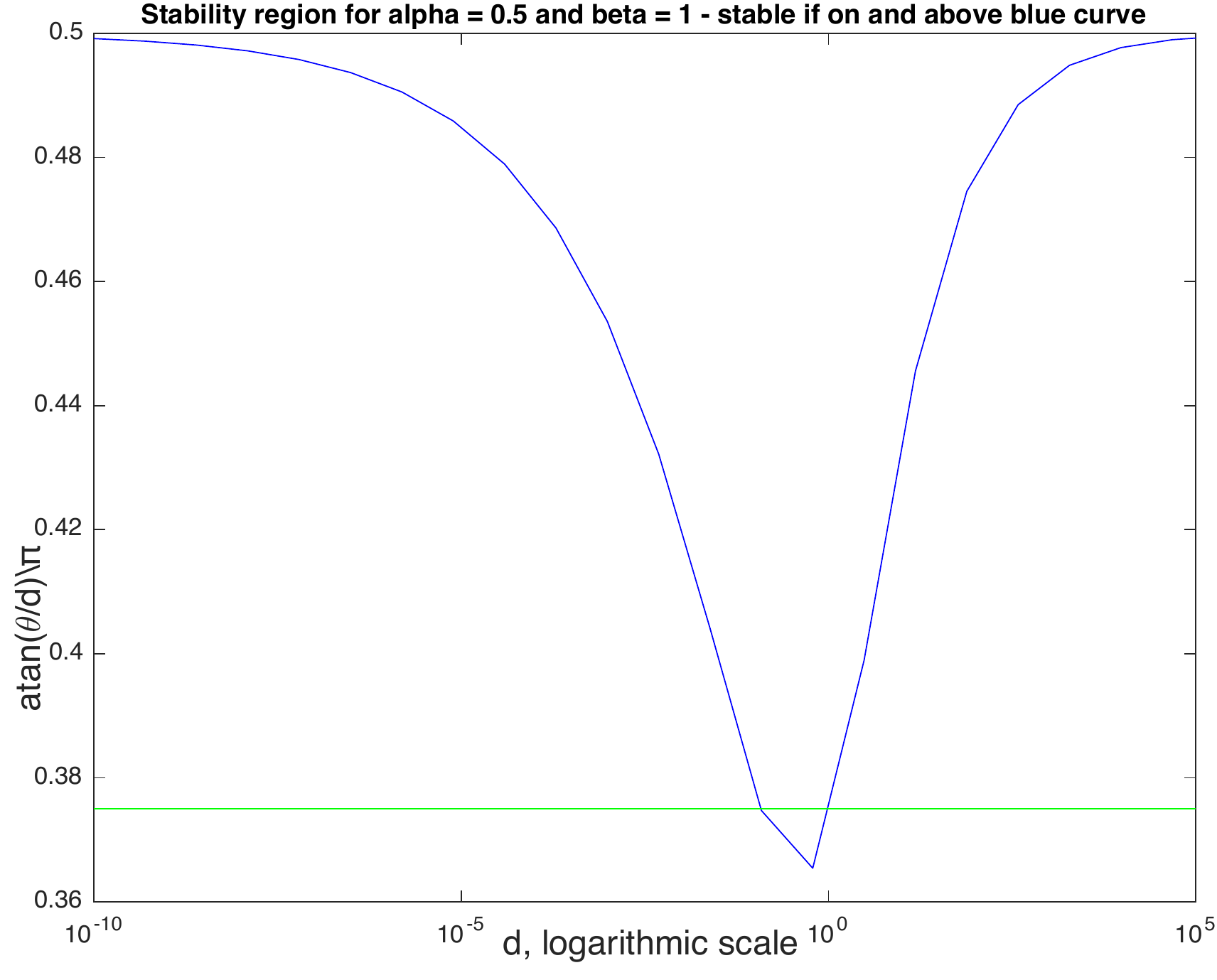} & \includegraphics[width=0.35\textwidth]{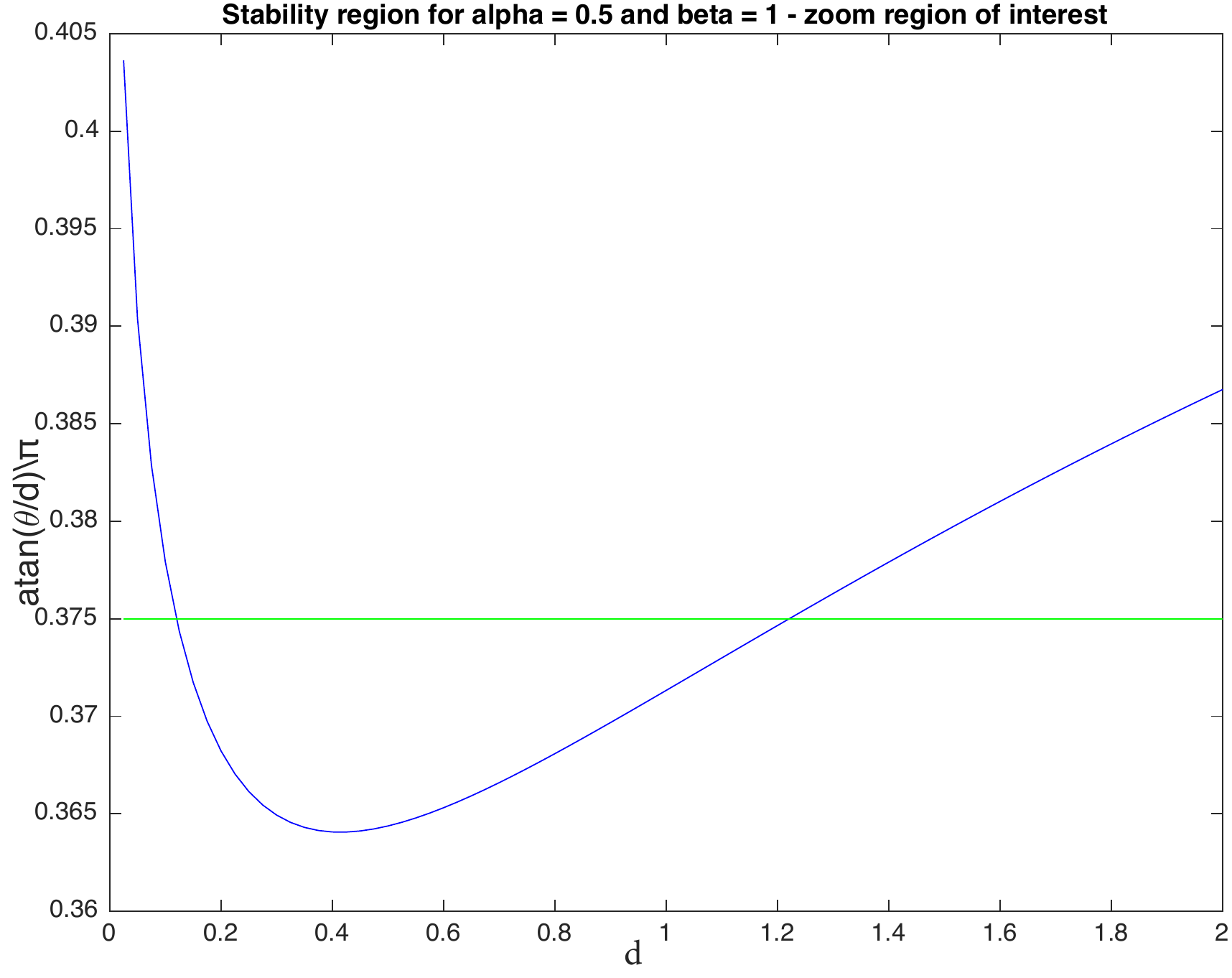}
\end{array}$$
\end{figure}

\begin{figure}[tbhp]
\captionsetup{width=0.75\textwidth}
\caption{Stability region, above the blue line, for choosing $d$ and $\theta$, when the eigenvalues of $A$ are $d \pm i \theta$, $\alpha = \frac{1}{3}, \beta = \frac{2}{3}.$ The logarithmic scale is explored in the right hand figure where the stability boundary dips below the angle $\frac{\pi}{4}$.} 
\label{fig:fig3}
\centering
$$\begin{array}{cc}
\includegraphics[width=0.35\textwidth]{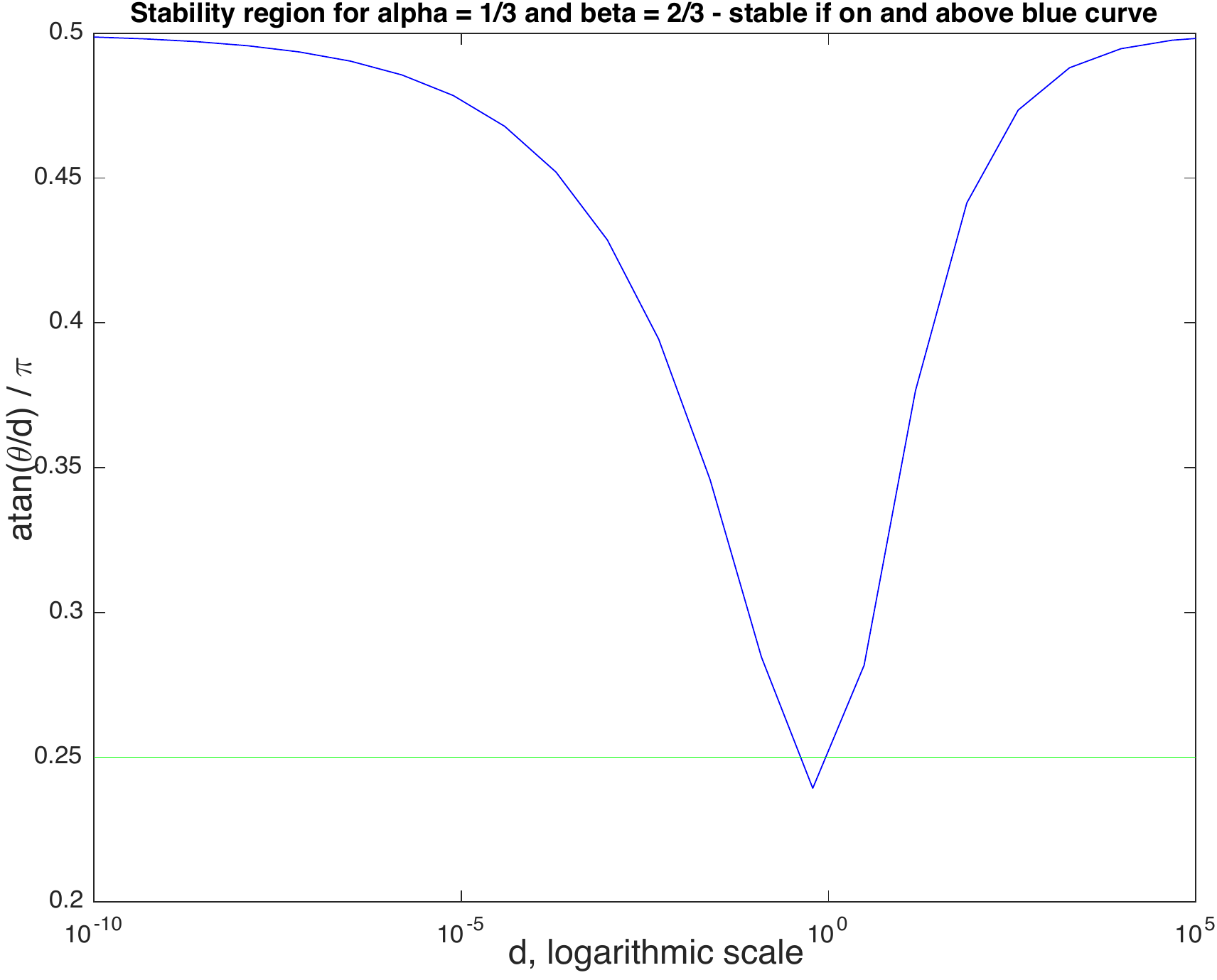} & \includegraphics[width=0.35\textwidth]{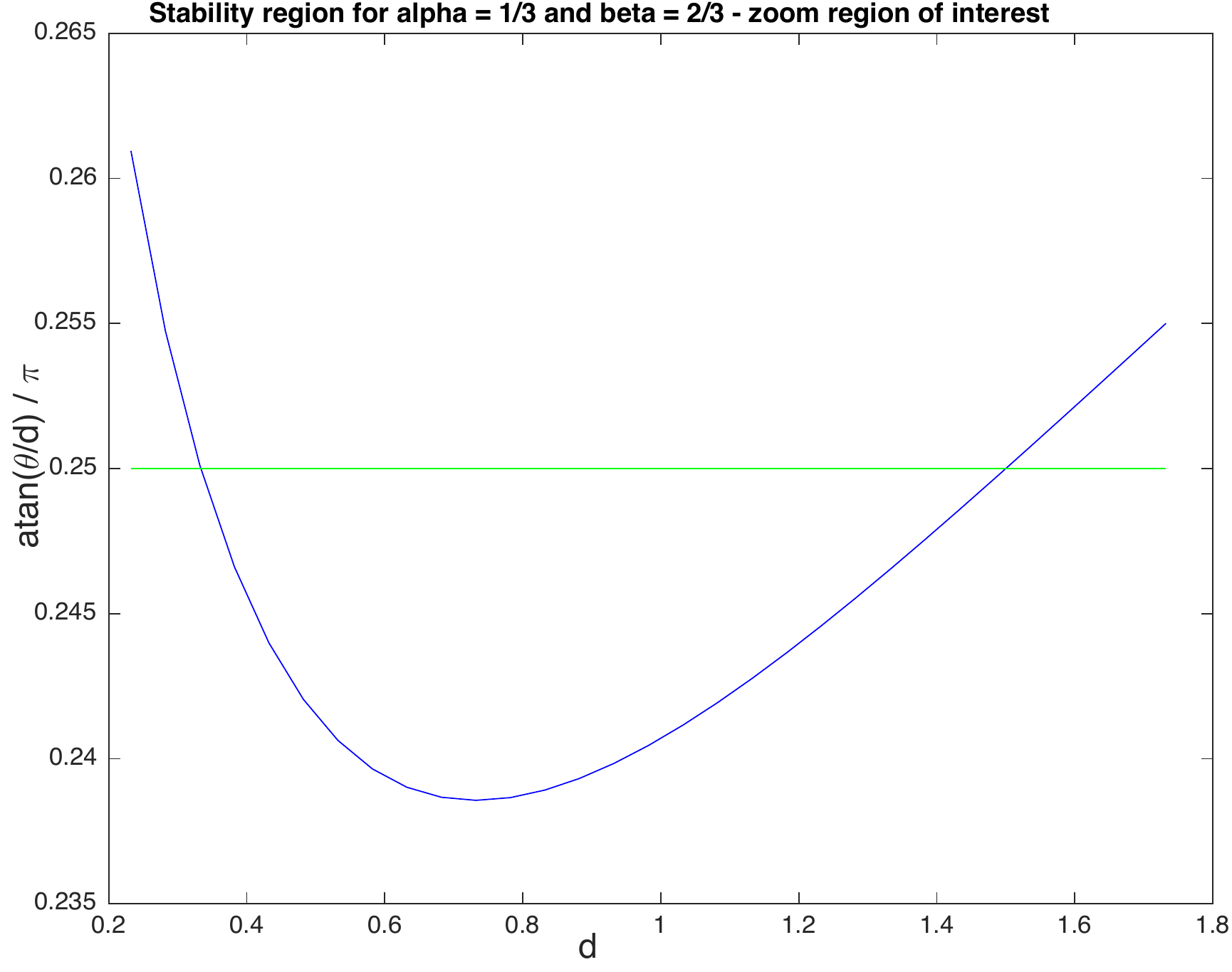}
\end{array}$$
\end{figure}

\begin{figure}[tbhp]
\captionsetup{width=0.75\textwidth}
\caption{System Dynamics with $(\alpha,\beta) = (\frac{1}{2},1)$, top, and $(\alpha,\beta) = (\frac{1}{3},\frac{2}{3})$, bottom. The left hand column shows sustained dynamics with $d = 1$ and $\theta$ chosen so that $(d,\theta)$ lies on the stability boundary. The right hand column corresponds to the same $d$ but 0.3 has been added to the $\theta$ value.} 
\label{fig:fig4}
\centering
$$\begin{array}{cc}
\includegraphics[width=0.35\textwidth]{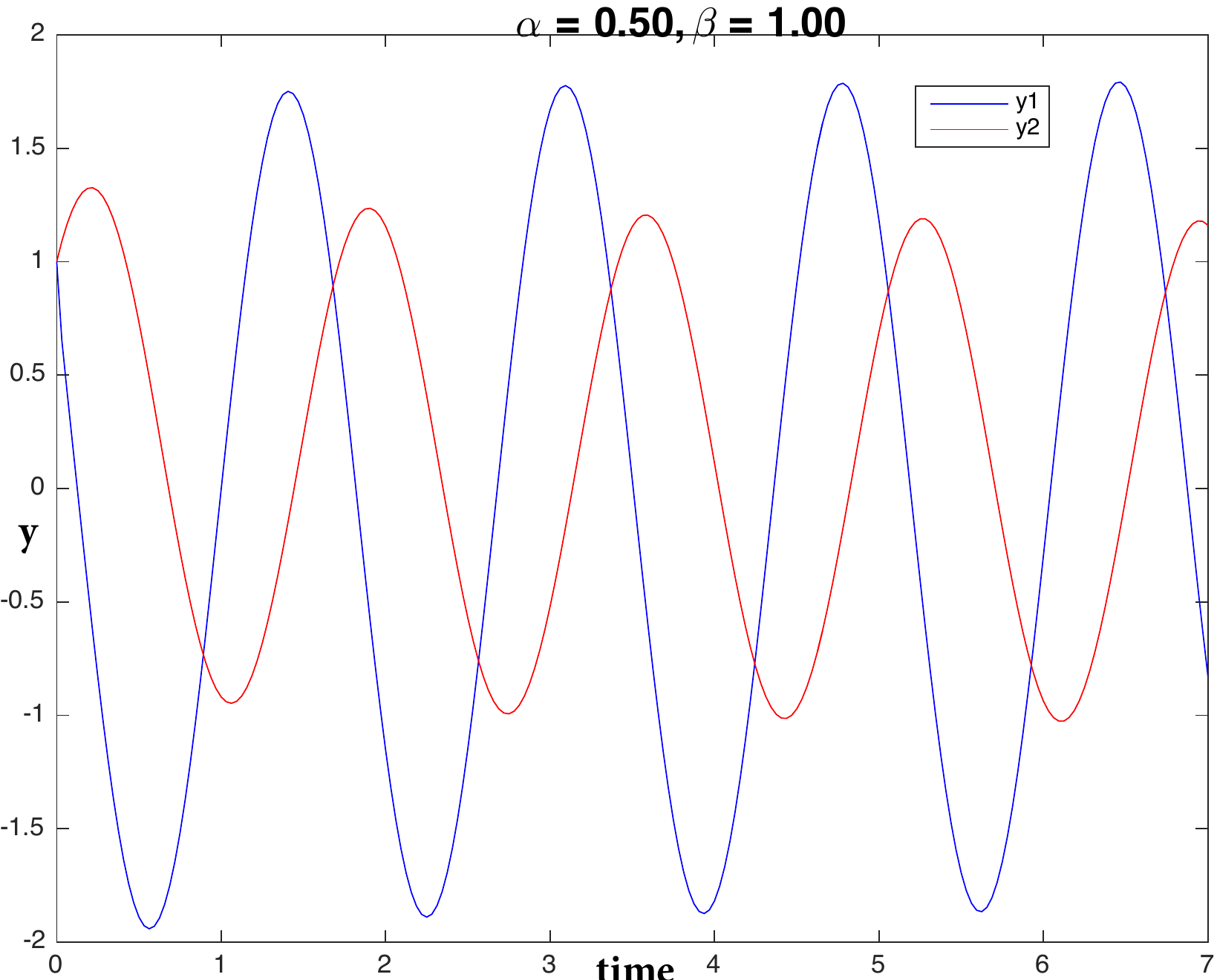} & \includegraphics[width=0.35\textwidth]{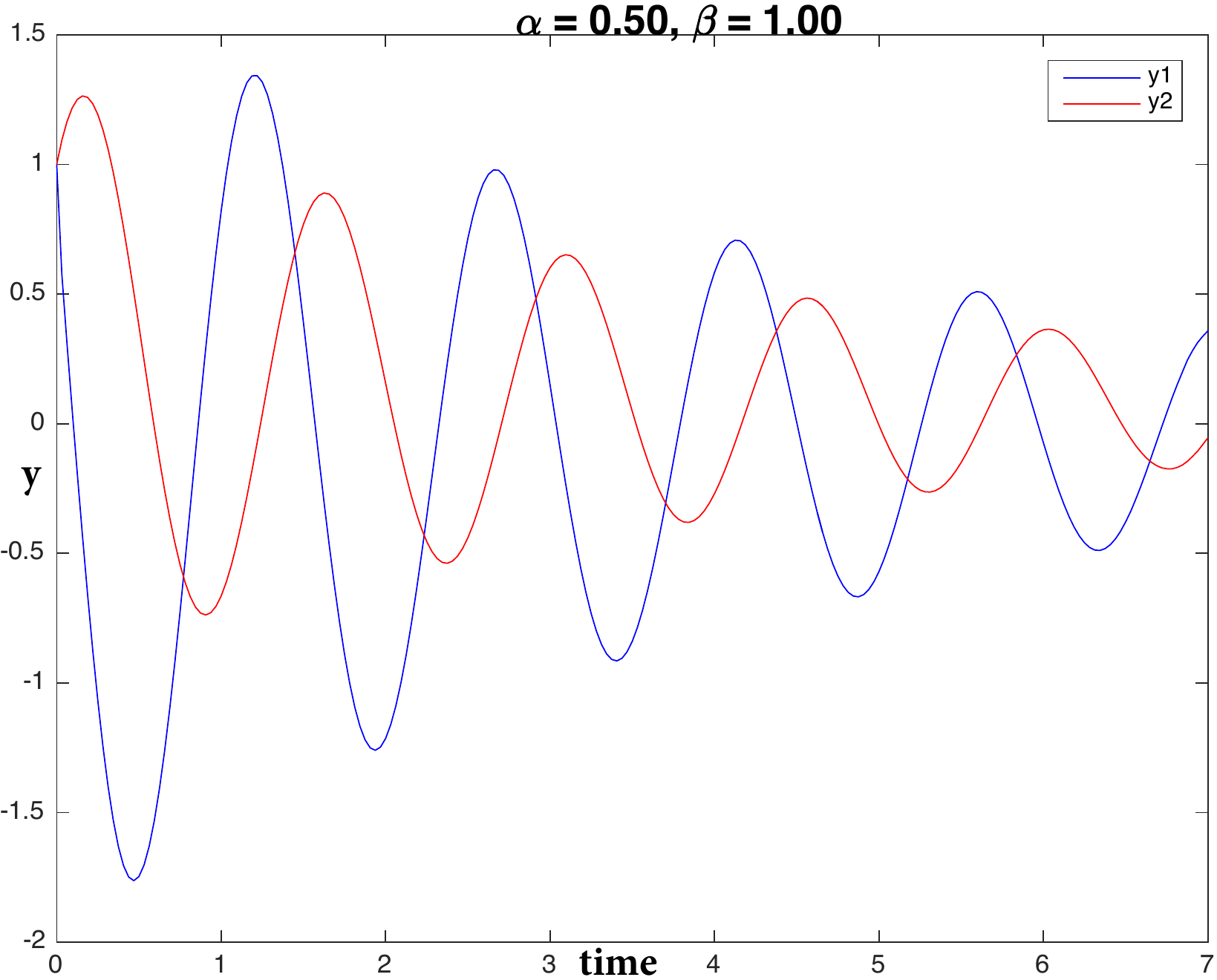} \\
& \\
\includegraphics[width=0.35\textwidth]{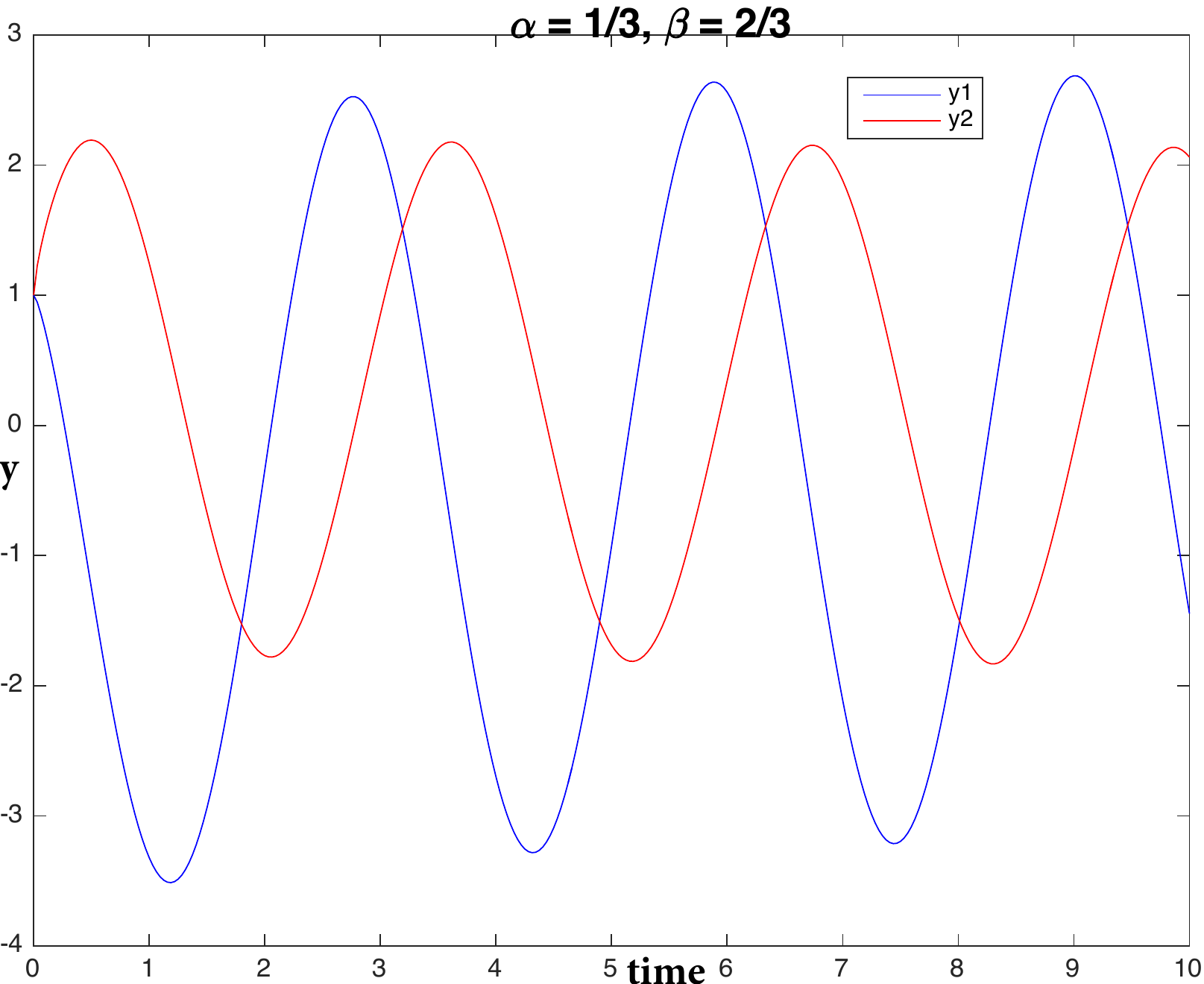} & \includegraphics[width=0.35\textwidth]{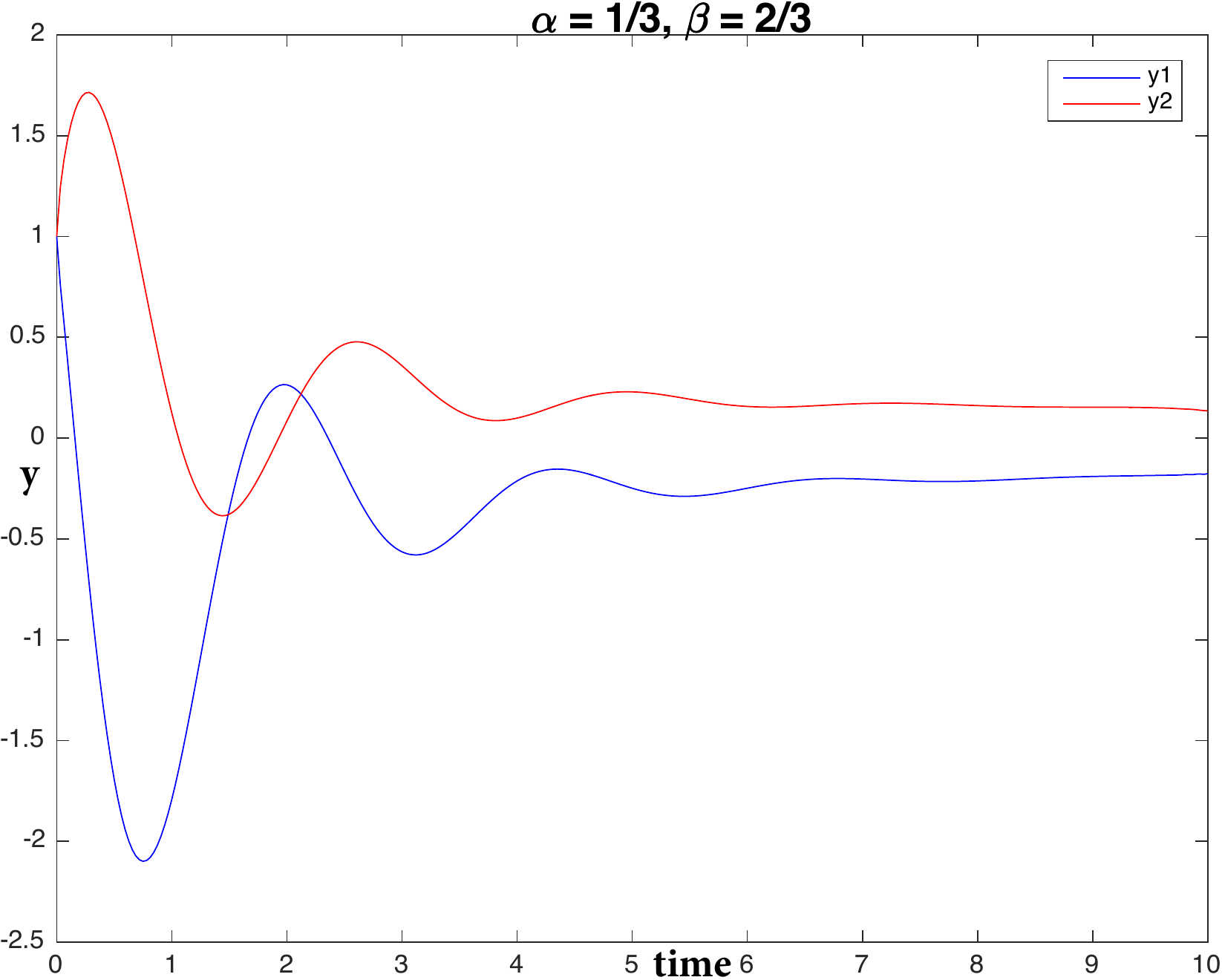} 
\end{array}$$
\end{figure}

\begin{figure}[tbhp]
\captionsetup{width=0.75\textwidth}
\caption{Phase Plots of $y_1$ versus $y_2$ for the decaying solutions in the right hand column of Figure \ref{fig:fig4}.} 
\label{fig:fig5}
\centering
$$\begin{array}{cc}
\includegraphics[width=0.35\textwidth]{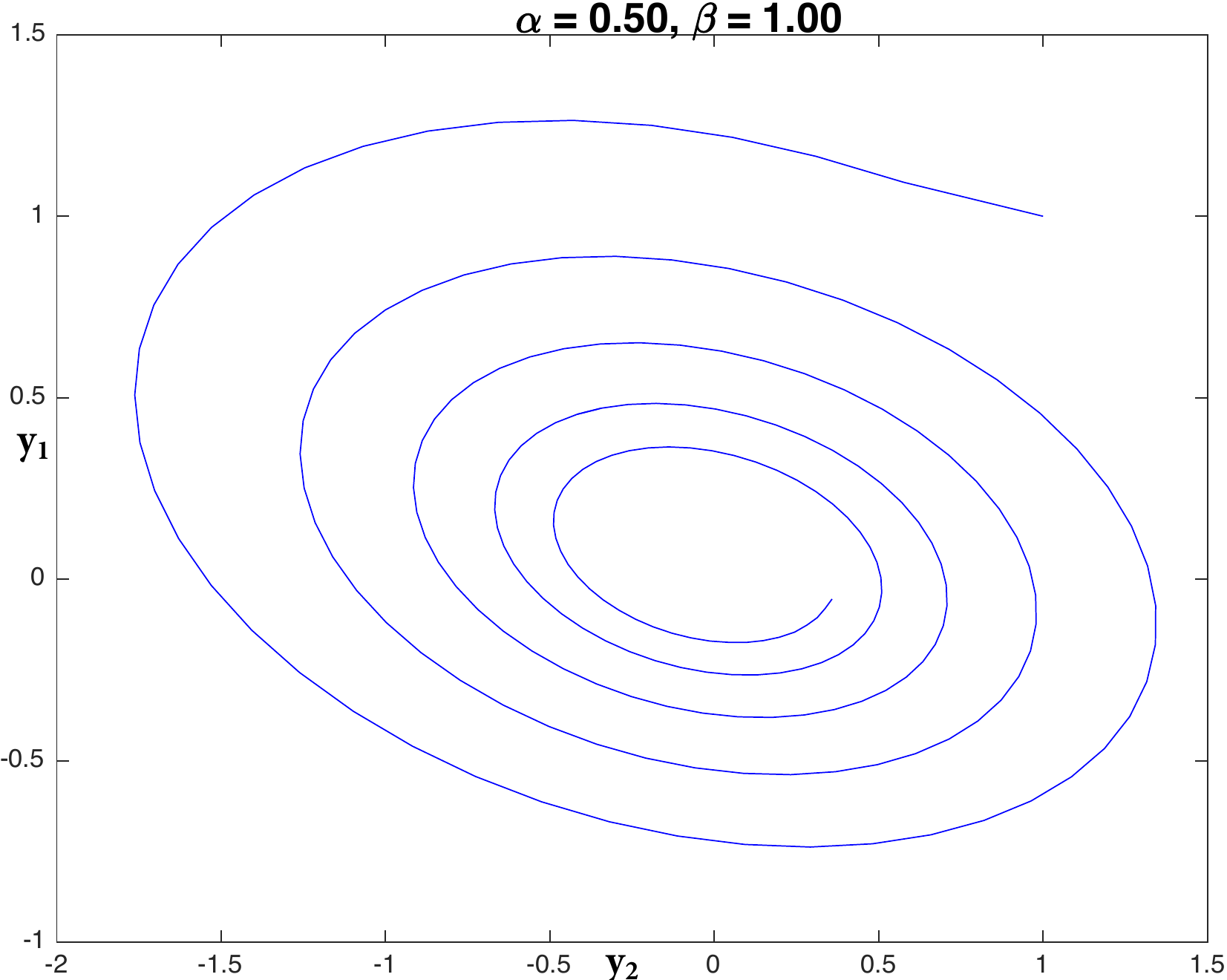} & \includegraphics[width=0.35\textwidth]{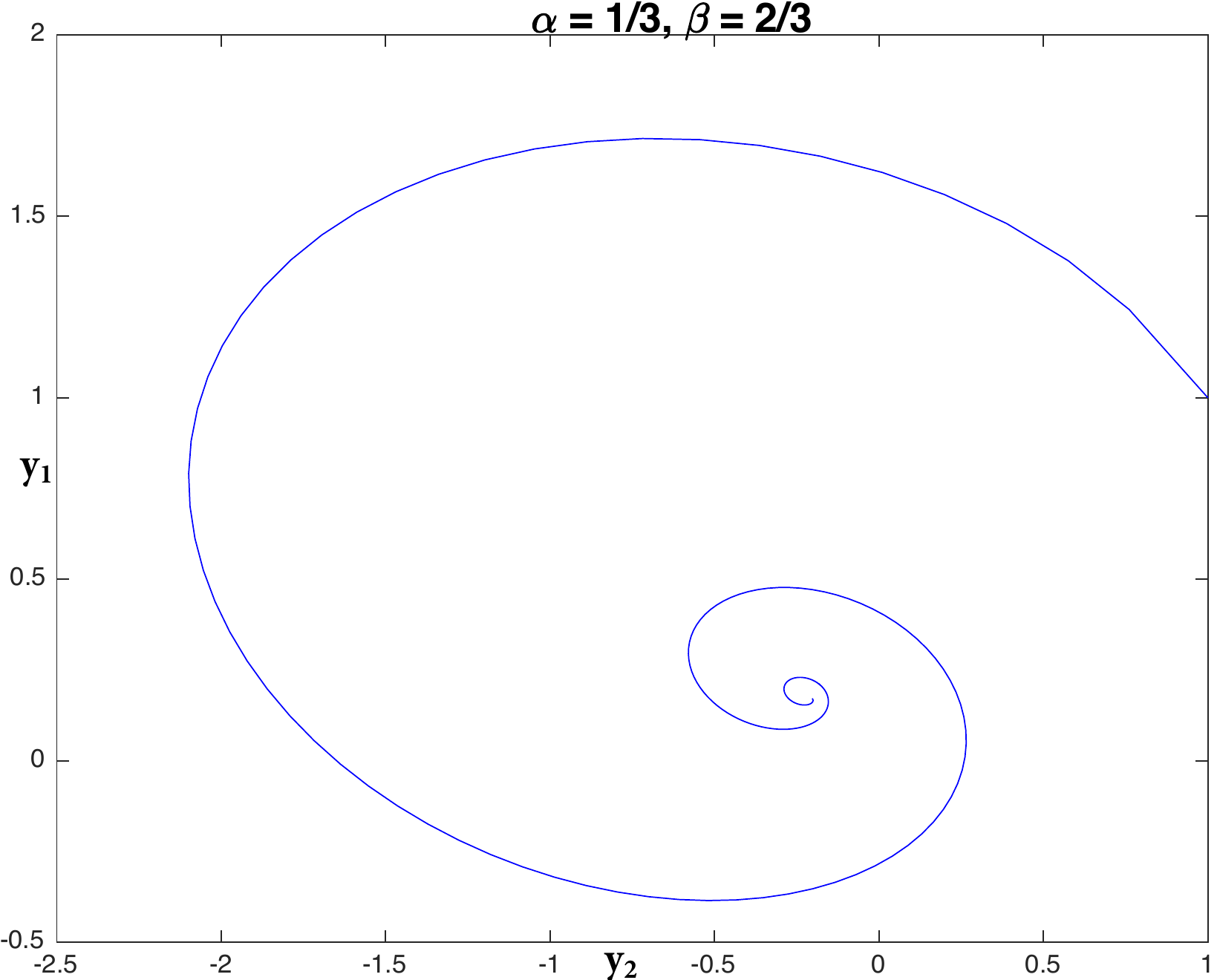} 
\end{array}$$
\end{figure}

\begin{figure}[htb]
\captionsetup{width=0.75\textwidth}
\caption{For $A$ given by (\ref{eq:diffA}) with $d=-1,\, \theta = \frac{1}{2}$ so that the eigenvalues are $-\frac{3}{2}, \, -\frac{1}{2}$, showing the effect of variation of $\alpha$ with fixed $\beta$ on the system dynamics.} 
\label{fig:fig6}
\centering
$$\begin{array}{cc}
\includegraphics[width=0.3\textwidth]{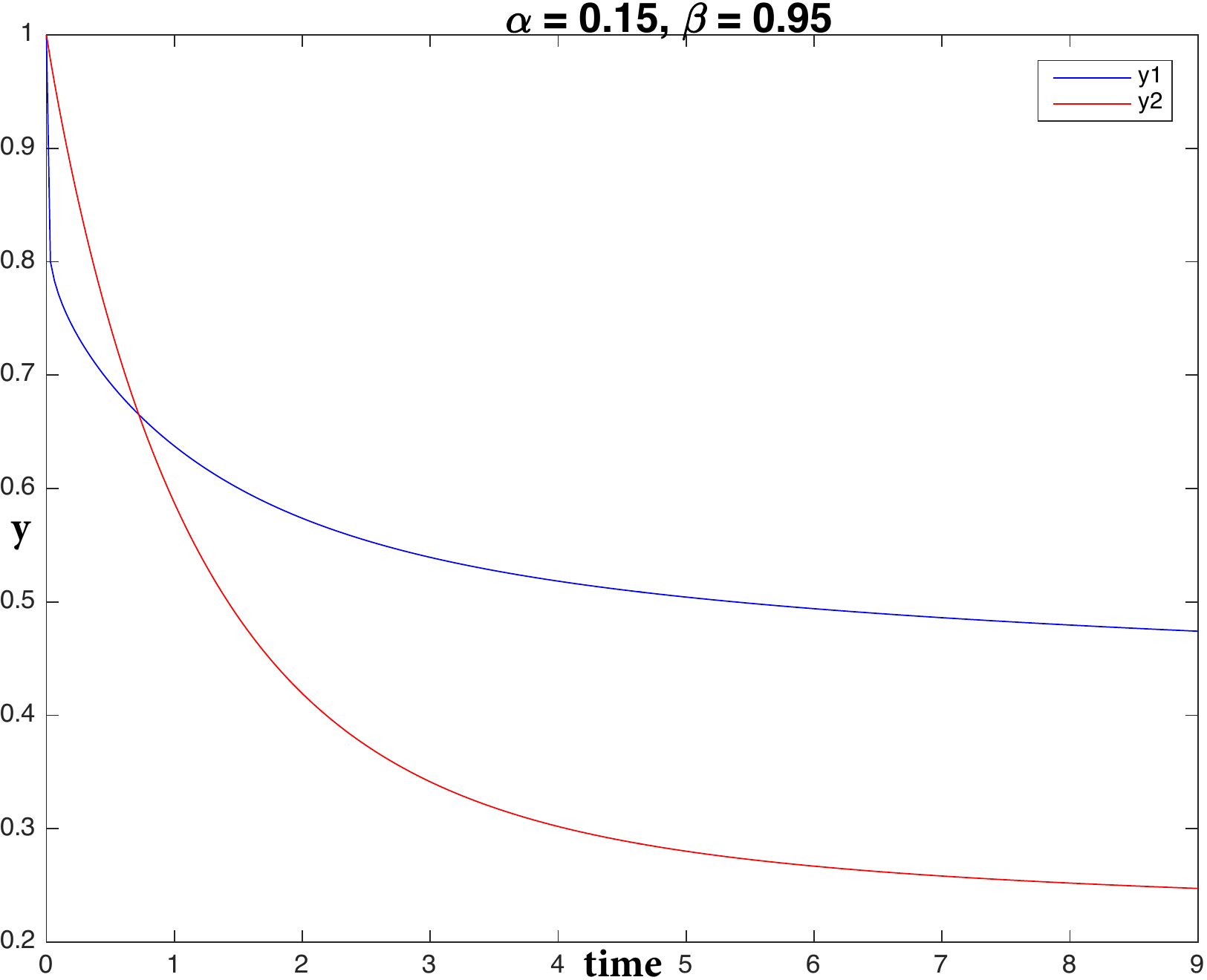} & \includegraphics[width=0.3\textwidth]{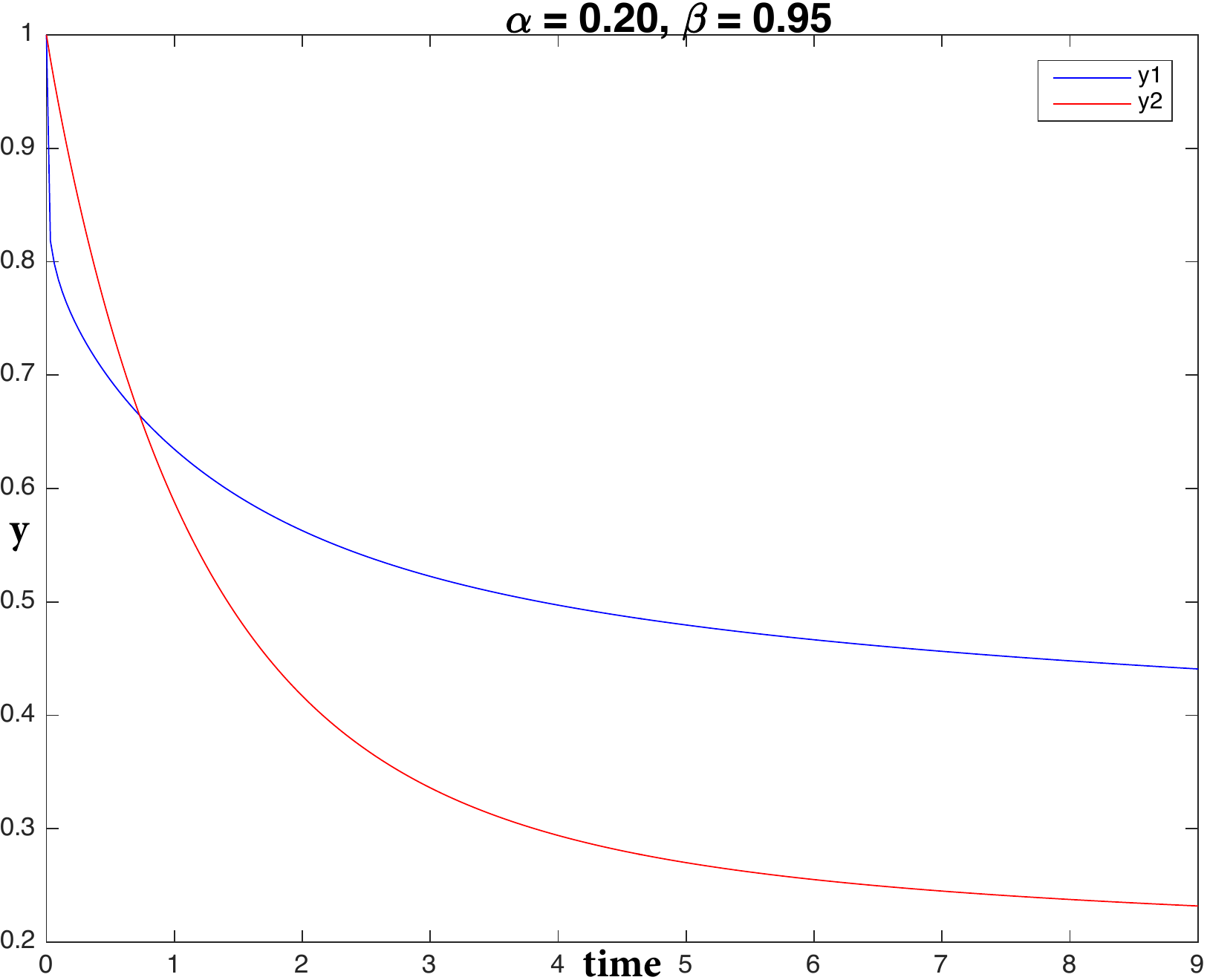} \\
& \\
\includegraphics[width=0.3\textwidth]{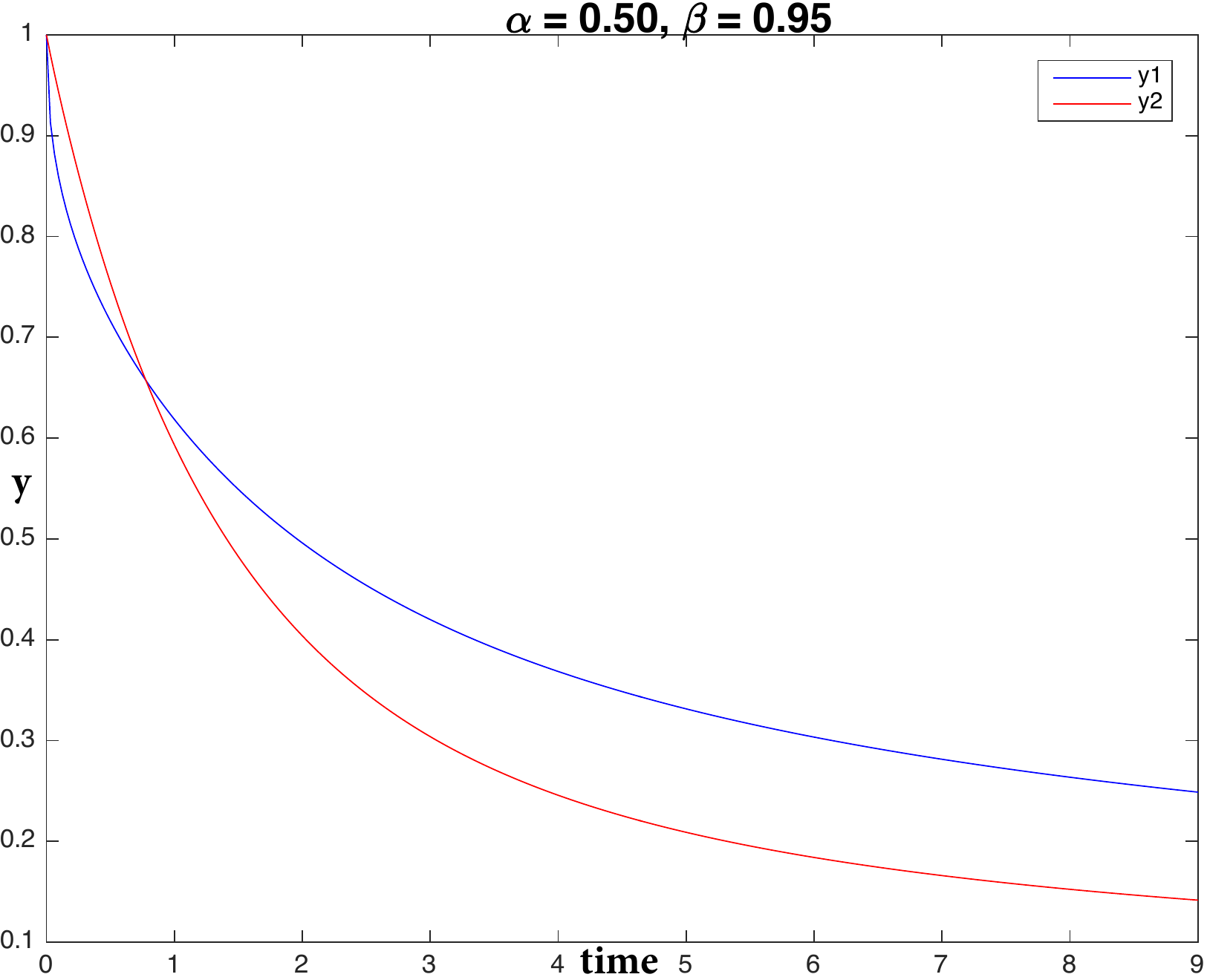} & \includegraphics[width=0.3\textwidth]{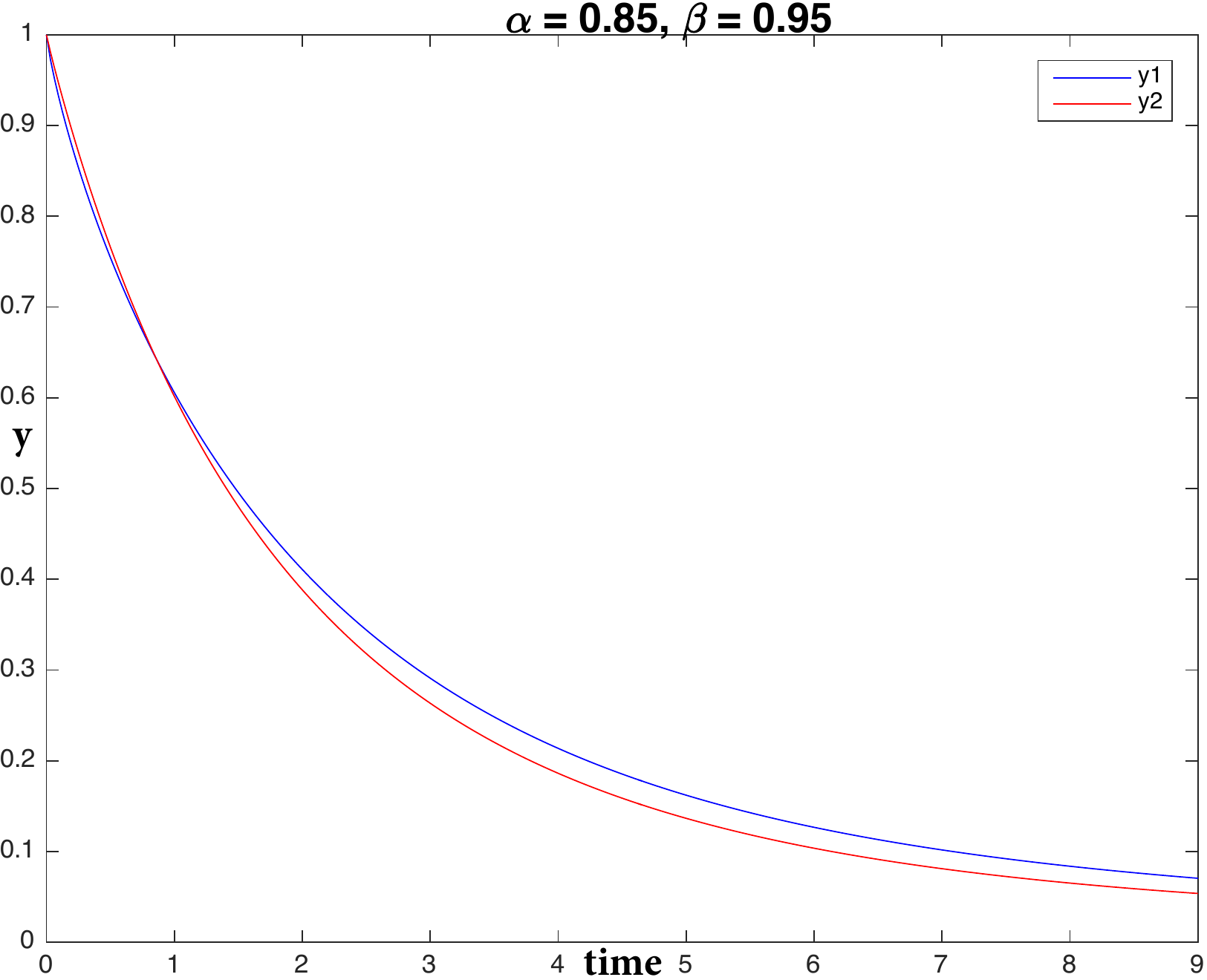} 
\end{array}$$
\end{figure}

In Figures \ref{fig:fig2} and \ref{fig:fig3} we plot the asymptotic stability boundary of the two dimensional, index-two problem given by (\ref{eq:1}) where
\begin{equation}
A = \left( \begin{array}{rr} d & -\theta \\ \theta & d  \end{array} \right), \quad d > 0
\label{eq:AMatrix}
\end{equation}
for the two cases considered in section 4, namely  $ \beta = 1, \> \alpha = \frac{1}{2}$ (Figure \ref{fig:fig2}) and $\beta=\frac{2}{3}, \> \alpha = \frac{1}{3}$ (Figure \ref{fig:fig3}).
Since the eigenvalues of $A$ are $d \pm i \theta$, we plot on the vertical axis the angle $\hat{\theta}$ in radians, where
$\hat{\theta} = \frac{1}{\pi} \arctan (\frac{\theta}{\lambda})$, as a function of $d$.
In Figure \ref{fig:fig2} we see that $\hat{\theta} \in (\frac{1}{4}, \frac{1}{2})$ corresponding to an angle lying between $45^\circ$ and $90^\circ$, as expected from the theory. We also plot the angle, in green, corresponding to the midpoint between these two extremes, i.e. $\frac{3}{8} \pi$. We see that for the most part the asymptotic stability angle lies above this midpoint except for the values of $d$, as shown in the right hand figure. 

In the case of Figure \ref{fig:fig3}, we give a similar plot as Figure \ref{fig:fig2}.  We also plot in green the midpoint between the two lines subtended by angles $\frac{1}{3}\pi$ and $\frac{1}{6}\pi$, namely $\frac{1}{4}\pi$. As with Figure \ref{fig:fig2} there is a small range of $d$ for which the asymptotic stability angle drops beneath $\frac{1}{4} \pi$. Furthermore, it is clear from Remarks 4(ii)  that as $\alpha$ and $\beta$ approach one another, the asymptotic stability boundary will be almost constant over increasingly longer periods of $d$ and will only asymptotically approach the angle $\frac{\pi}{2}$ for very small and very large values of $d$ - see Remark part (ii).

In Figure \ref{fig:fig4} we confirm the asymptotic stability analysis showing sustained and decaying oscillations with $\alpha = \frac{1}{2}, \beta = 1$ (top panel) and $\alpha = \frac{1}{3}, \beta = \frac{2}{3}$ (bottom panel). In all four cases, $d = 1$ while for the top panel we take $\theta = \sqrt{2(1+\sqrt{3})}, \> \theta = \sqrt{2(1+\sqrt{3})}+0.3$, while for the bottom panel we take $\theta = \frac{\sqrt{3}}{4}\sqrt{ \sqrt{33}-1}, \> \theta = \frac{\sqrt{3}}{4}\sqrt{ \sqrt{33}-1} + 0.3$.

In Figure \ref{fig:fig5} we present phase plots of $y_1$ versus $y_2$ for the two decaying oscillations cases. The figures confirm our theoretical results on the asymptotic stabiity boundary and also show the effects that the fractional indices have on the period of the solutions. As $\alpha$ approaches $\beta$ we expect the oscillatory behaviour to disappear.

Finally, in Figure \ref{fig:fig6} we consider the problem 
\begin{equation}
A = \left( \begin{array}{cc} d & \theta \\ \theta & d \end{array} \right), \quad d < 0
\label{eq:diffA}
\end{equation}
in which case the eigenvalues of $A$ are $d \pm \theta$. We take $d = -1, \, \theta = \frac{1}{2}$ and present the solutions for four pairs of indices, namely $(\alpha, \beta)  = (0.85, 0.95), \, (0.5, 0.95), \, (0.2, 0.05)$, $(0.15, 0.95)$. The simulations show that the components of the solution $y_1$ and $y_2$ seem to pick up ``energy" from one another due to the coupling and that as the distance between $\alpha$ and $\beta$ grows there is a greater separation between the two components. Finally, as $\alpha$ gets smaller, the solutions appear to ``flat-line" more quickly.



\section{Conclusions}

In this paper we have studied mixed index fractional differential equations with coupling between the different components. We find an analytical expression for the solution of the linear system that generalises the Mittag-Leffler expansion of a matrix and the solution of linear sequential fractional differential equations. We can use this result to derive new numerical methods that generalise the concept of exponential methods used in the approximation of the Mittag-Leffler matrix function, see \cite{ref20, ref21, ref21a},  for example, and exponential integrators \cite{ref22}, \cite{ref23}. The second element would deal with developing numerical techniques for the integration component that incorporates the integral of a function times a Green function. We also use Laplace transform techniques to find the asymptotic stability domain in terms of the eigenvalues of the defining linear system. Finally we have also used Laplace transforms to get  analytical expansions of the mixed index problem in terms of a sum of Mittag-Leffler or generalised Mittag-Leffler functions, in the case that the fractional indices are rational.

\section{Acknowledgements}
We would like to thank  Dr Alfonso Bueno-Orovio in the Department of Computer Science, University of Oxford, for many discussions about fractional differential equations.

\end{document}